\documentclass{amsart}

\usepackage{delimset, amssymb, enumitem, mathtools, keytheorems}
\setlist{itemsep=4pt, topsep=0pt, leftmargin=17pt, listparindent=11pt}

\usepackage{xcolor}
\definecolor{BIT}{cmyk}{1, 0, 1, 0}
\usepackage[colorlinks, citecolor=BIT, pagebackref, ocgcolorlinks]{hyperref}

\usepackage[capitalize]{cleveref}
\crefname{def}{Def.}{Defs.}
\Crefname{def}{Definition}{Definitions}
\creflabelformat{def}{~\upshape(#2#1#3)}
\crefname{ineq}{Ineq.}{Ineqs.}
\Crefname{ineq}{Inequality}{Inequalities}
\creflabelformat{ineq}{~\upshape(#2#1#3)}
\crefname{conjecture}{Conjecture}{Conjectures}

\AddToHook{env/conjecture/begin}{\crefalias{theorem}{conjecture}}
\AddToHook{env/lemma/begin}{\crefalias{theorem}{lemma}}
\AddToHook{env/proposition/begin}{\crefalias{theorem}{proposition}}

\makeatletter
\newcommand\creflabel[2][\@currentcounter]{%
 \crefalias{\@currentcounter}{#1}\label{#2}}
\makeatother

\newtheorem{theorem}{Theorem}[section]
\newkeytheorem{corollary}[name=Corollary, sibling=theorem]
\newkeytheorem{conjecture}[name=Conjecture, sibling=theorem]
\newkeytheorem{lemma}[name=Lemma, sibling=theorem]
\newkeytheorem{proposition}[name=Proposition, sibling=theorem]

\numberwithin{equation}{section}

\allowbreak
\allowdisplaybreaks

\usepackage[numbers, sort&compress, nonamebreak, merge, elide, longnamesfirst]{natbib}

\usepackage{url}

\usepackage{tikz}
\usetikzlibrary{decorations.pathreplacing, calc, positioning}
\tikzset{
edge/.style={semithick},
ball/.style={shape=circle, minimum size=1mm, ball color=black, inner sep=0.5},
ellipsis/.style={shape=circle, fill, inner sep=.5}}
\def\r{1}
\def\eps{.2}

\usepackage{xcolor}

\title{Chromatic symmetric functions of conjoined graphs}

\author[E.Y.J. Qi]{Ethan Y.J. Qi}
\address[Ethan Y.J. Qi]{School of Mathematics and Statistics, Beijing Institute of Technology, Beijing 102400, P. R. China.}
\email{ethane@bit.edu.cn}

\author[D.Q.B.~Tang]{Davion Q.B. Tang}
\address[Davion Q.B. Tang]{School of Mathematics and Statistics, Beijing Institute of Technology, Beijing 102400, P. R. China}
\email{davion@bit.edu.cn}

\author[D.G.L.~Wang]{David G.L. Wang$^*$}
\address[David G.L. Wang]{School of Mathematics and Statistics \& MIIT Key Laboratory of Mathematical Theory and Computation in Information Security, Beijing Institute of Technology, Beijing 102400, P. R. China.}
\email{glw@bit.edu.cn}
\thanks{$^*$Wang is the corresponding author, and is supported by the NSFC (Grant No.\ 12171034).}


\begin{document}

\bibliographystyle{abbrvnat}

\begin{abstract}
We introduce path-conjoined graphs defined for two rooted graphs by joining their roots with a path, and investigate the chromatic symmetric functions of its two generalizations: spider-conjoined graphs and chain-conjoined graphs. By using the composition method developed by Zhou and the third author recently, we obtain neat positive $e_I$-expansions for the chromatic symmetric functions of clique-path-cycle graphs, path-clique-path graphs, and clique-clique-path graphs. We pose the $e$-positivity conjecture for hat-chains.
\end{abstract}
\maketitle
%\tableofcontents

\section{Introduction}
The chromatic symmetric function of a graph was introduced by \citet{Sta95} as a generalization of Birkhoff's chromatic polynomial.
It attracted a considerable number of studies 
from various branches of mathematics, 
including combinatorics, representation and knot theory, and from mathematical physics, see \citet{NW99,Loe07,CKL20,Kato24X}.

One of the most popular topics in this field 
is \citeauthor{SS93}'s conjecture~\cite{SS93}.
In virtue of \citeauthor{Gua13X}'s reduction \cite{Gua13X},
it is equivalent to the $e$-positivity of 
unit interval graphs.
More precisely,
the chromatic symmetric function of any unit interval graph
has only nonnegative $e$-coefficients.
The conjecture charmed representation experts
because of a close relation between
chromatic quasisymmetric functions that were introduced by \citet{SW16}
and the representations of symmetric groups
on cohomology of regular semisimple Hessenberg varieties,
see \citet{BC18}.
In October 2024,  
\citet{Hik25} resolved \citeauthor{SS93}'s conjecture.

A natural generalization of the conjecture
is to characterize all $e$-positive graphs, as mentioned by \citet[Section 5]{Sta95}.
Though seemingly unachievable is it, 
more and more families of graphs have been confirmed to be $e$-positive.
The most noted $e$-positive graphs are complete graphs, paths, cycles,
and graphs whose independent numbers are $2$.
A bit less well-known $e$-positive graphs are lollipops,
tadpoles and generalized bulls, 
which can be derived from a direct application of 
\citeauthor{GS01}'s 
$(e)$-positivity strengthening, see \cite[Theorem 7.6]{GS01},
and also \citet{Dv18,LLWY21}.
Here strengthening means that the $(e)$-positivity is stronger than the $e$-positivity.
More $e$-positive graphs can be found from \citet{Tsu18,CH19,FHM19,HHT19,LY21,WW23-JAC,Wang25,BCCCGKKLLS25}.

Along the way of confirming the $e$-positivity of graphs arose various approaches.
Notable routes include 
(i) \citeauthor{Sta95}'s generating function and recurrence relation method,
(ii) the strategy of expressing a chromatic symmetric function 
as a sum of $e$-positive symmetric functions, 
often as a sum of known $e$-positive chromatic symmetric functions; 
a typical example is the proof for the $3$-spiders of the form $S(n+1,n,1)$, 
see \citet[Page 2681]{DFv20}, 
(iii) \citeauthor{GS01}'s $(e)$-positivity strengthening 
that takes advantage of the algebra $\mathrm{NCSym}$ 
of symmetric functions in noncommutative variables, 
(iv) understanding the $e$-coefficients of a chromatic symmetric function
combinatorially
and constructing sign-reversing involutions, see \citet{Tom25-CT} for instance,
(v) the composition method which is inspired from \citet{SW16}'s formula 
for the chromatic symmetric function of paths, initiated implicitly 
by \citet{TW23X} in proving the Schur positivity 
of the $3$-spiders of the forms $S(a,2,1)$ and $S(a,4,1)$, 
and formulated by \citet{WZ25}.

In this paper,
we continue explore the power
of the composition method.
A key idea of the composition method is
to extend the defining domain of the elementary symmetric functions $e_\lambda$
from partitions $\lambda$
to compositions $I$. 
An \emph{$e_I$-expansion} of a symmetric function $f$ is an expression
\begin{equation}\label{def:eI-expansion}
f=\sum_{I\vDash n} c_I e_I
\end{equation}
such that 
\begin{equation}\label{def2:eI-expansion}
f=\sum_{\lambda\vdash n}
\brk4{\sum_{\rho(I)=\lambda}c_I}
e_\lambda,
\end{equation}
where $\rho(I)$ is the partition obtained from the composition $I$ by rearranging its parts.
We call \cref{def:eI-expansion} a \emph{positive $e_I$-expansion}
if $c_I\ge 0$ for all compositions $I$ of $n$.
In order to show that~$f$ is $e$-positive,
it suffices to find a positive $e_I$-expansion of $f$.
Every $e$-positive symmetric function has definitely a positive $e_I$-expansion,
while the non-uniqueness of $e_I$-expansions suggests
a flexibility in invoking combinatorial ideas to expose 
or construct a positive $e_I$-expansion. 

One family of the foregoing $e$-positive graphs that are unit interval graphs is \emph{$K$-chains}, see \cref{sec:pre}.
Its $e$-positivity was due to \citet{GS01},
and rediscovered by \citet{Tom24,MPW24X}. 
We would like to seek for a neat positive $e_I$-expansion for $K$-chains.
Here the word neat means that
essentially less compositions appear in the expansion,
or aesthetically
less technical conditions are required for the compositions 
that appear.
This has been done
for 
two particular kinds of $K$-chains: barbells and generalized bulls, see \cite{WZ25}.
In this paper, 
we continue this course and light on another two $K$-chains: path-clique-paths (PKPs)
and clique-clique-paths (KKPs), see \cref{thm:PKP,thm:KKP} respectively.
In fact, all generalized bulls are PKPs.
In the same fashion, 
we find a positive $e_I$-expansions for clique-path-cycles (KPCs), 
see \cref{thm:KPC}.
This is the first positive $e_I$-expansion of KPC graphs
to the best of our knowledge.
It is transparent that 
the composition method is distinct
from the involution one,
and our positive $e_I$-expansions for PKPs and KKPs are 
different from \citeauthor{Tom24}'s, see~\cref{thm:K-chain}.
We should mention that the $e$-positivity of all these graphs
are derivable by using \citeauthor{GS01}'s $(e)$-positivity routine.

Ad hoc algebraic combinatorial techniques are called 
for in coping with PKPs, KKPs and KPCs,
nevertheless, these graphs have a common structure:
a path part as a bridge that links some \emph{node graphs}
such like cliques or cycles.
We are therefore led to the graph operation of \emph{path-conjoining},
see \cref{sec:path-conjoined}. 
We then generalize it
to spider- and chain-conjoining operations,
and express the chromatic symmetric function of these amalgamated graphs
in terms of those of the node graphs.

This paper is organized as follows.
\cref{sec:pre}
consists of preliminary results that we will need in the sequel.
In \cref{sec:path-conjoined},
we study the chromatic symmetric function of path-conjoined graphs, 
and specify our results to those with a clique node and those with a cycle node. 
\cref{sec:spider-conjoined} is devoted 
to spider-conjoined graphs. In particular,
we investigate some pineapple graphs
which appeared in \citet{TSH16}, see \cref{cor:pineapple:mab}.
In \cref{sec:chain-conjoined},
we consider $2$-chain-conjoined graphs with the middle node being
a clique or a cycle, see \cref{thm:KGH,thm:cycle2path:AB}.
We end this paper with the $e$-positivity conjecture for hat-chain graphs.

\section{Preliminaries}\label{sec:pre}

This section contains 
basic knowledge on chromatic symmetric functions.
We adopt notion and notations of \citet{Sta11B}.
Let $n$ be a positive integer. 
A \emph{composition} of $n$ is 
a sequence of positive integers with sum~$n$,
commonly denoted 
$I=i_1 \dotsm i_l\vDash n$, 
with \emph{size} $\abs{I}=n$, \emph{length} $\ell(I)=l$, 
and \emph{parts} $i_1,\dots,i_l$.
A \emph{prefix} of $I$
is a composition of the form $i_1\dotsm i_k$, where $0\le k\le \ell(I)$.
We write $I=v^s$ if all parts have the same value $v$.
For notational convenience, 
we denote the $k$th last part $i_{l+1-k}$
by~$i_{-k}$,
and denote the composition obtained 
by removing the $k$th part by~$I\backslash i_k$, i.e.,
\[
I\backslash i_k
=i_1\dotsm i_{k-1}i_{k+1}\dotsm i_{-1}.
\]
When a capital letter like $I$ and $J$
stands for a composition,
we use its small letter counterpart
with integer subscripts to denote the parts.
A \emph{partition}\index{partition} of $n$
is a multiset of positive integers $\lambda_i$ with sum~$n$,
denoted 
$\lambda=\lambda_1\lambda_2\dotsm
\vdash n$,
where $\lambda_1\ge \lambda_2\ge\dotsm\ge 1$.

A \emph{symmetric function} of homogeneous degree $n$ over
the field $\mathbb{Q}$ of rational numbers is a formal power series
\[
f(x_1, x_2, \dots)
=\sum_{
\lambda=\lambda_1 \lambda_2 \dotsm \vdash n}
c_\lambda \cdotp
x_1^{\lambda_1} 
x_2^{\lambda_2} 
\dotsm, 
\]
where $c_\lambda \in \mathbb Q$,
such that $f\left(x_1, x_2, \ldots\right)=f\left(x_{\pi(1)}, x_{\pi(2)}, \ldots\right)$
for any permutation $\pi$.
One basis of 
the vector space of homogeneous symmetric functions of degree $n$ over $\mathbb{Q}$
is the set of elementary symmetric functions 
$e_\lambda
=e_{\lambda_1}
e_{\lambda_2}\dotsm$, where
\[
e_k
=\sum_{1\le i_1<\dots<i_k} 
x_{i_1} \dotsm x_{i_k}.
\]
A symmetric function $f$ is said to be \emph{$e$-positive}
if every $e_\lambda$-coefficient is nonnegative.

For any composition $I$, there is a unique partition $\rho(I)$
that consists of the parts of $I$.
As a consequence, 
one may extend the defining domain of elementary symmetric functions 
from partitions to compositions as
$e_I=e_{\rho(I)}$.
An \emph{$e_I$-expansion} of a symmetric function $f$ is an expression
\cref{def:eI-expansion} 
that satisfies \cref{def2:eI-expansion}.
It is said to be \emph{positive} 
if $c_I\ge 0$ for all compositions $I\vDash n$.
In order to show that~$f$ is $e$-positive,
it suffices to find a positive $e_I$-expansion of $f$.

\citet{Sta95} introduced the chromatic symmetric function for a simple graph $G$ as
\[
X_G
=\sum_\kappa \prod_{v \in V(G)} x_{\kappa(v)},
\]
where $\kappa$ runs over proper colorings of~$G$. 
It is a homogeneous symmetric function of degree $\abs{G}$,
the order of $G$.
For instance, the chromatic symmetric function of the complete graph~$K_n$ is $X_{K_n}=n!e_n$.

\Citet[Theorem 3.1, Corollaries 3.2 and 3.3]{OS14} established the \emph{triple-deletion} properties for chromatic symmetric functions as follows.

\begin{proposition}[\citeauthor{OS14}, the triple-deletion properties]\label{prop:3del}
Let $G$ be a graph with a stable set $T$ of order $3$.
Denote the edges linking the vertices in $T$ by $e_1$, $e_2$ and $e_3$.
For any set $S\subseteq \{1,2,3\}$, 
denote by $G_S$ the graph with vertex set~$V(G)$
and edge set $E(G)\cup\{e_j\colon j\in S\}$.
Then 
\begin{align}
\label{fml:3del:12}
X_{G_{12}}
&=X_{G_1}+X_{G_{23}}-X_{G_3},
\quad\text{and}\\
\label{fml:3del:123}
X_{G_{123}}
&=X_{G_{13}}+X_{G_{23}}-X_{G_3}.
\end{align}
\end{proposition}

Using \cref{prop:3del},
\citet[Proposition 6]{AWv24}
showed that certain sequence of chromatic symmetric functions 
forms an arithmetic progression.
A \emph{clique} of a graph~$G$
is a subset of vertices of $G$ 
such that every two distinct vertices in the clique are adjacent.
The \emph{closed neighborhood} of a vertex $v$ of $G$
is the union 
of the set of vertices that are adjacent to $v$ and~$v$ itself.

\begin{proposition}[\citeauthor{AWv24}, the arithmetic progression property]\label{prop:AP}
Suppose that $G$ contains a clique on the vertices $v_1, \dots, v_k$, 
which have the same closed neighborhoods.
Let $w$ be a vertex of~$G$ that is not adjacent to any vertex $v_j$.
Then the sequence $X_{G_0}, \dots, X_{G_k}$ forms an arithmetic progression, 
where $G_0=G$ and $G_j$ for $1\le j\le k$ is obtained from~$G_{j-1}$ by adding the edge $w v_j$.
\end{proposition}

In applying \cref{prop:AP}, 
we define
a \emph{tribe-vertex pair} of a graph $G$ to be a pair $(K,x)$
in which
\begin{enumerate}
\item
$K$ is a clique of $G$,
\item
$x$ is a vertex that may or may not be in $G$, and
\item
all vertices in $K$ 
have the same closed neighborhoods in the graph $G-x$.
\end{enumerate}

\begin{corollary}\label{cor:AP:add}
Let $(K,x)$ be a tribe-vertex pair of a graph $G$.
Let~$G_i$ be the graph obtained from~$G$ by
connecting exactly $i$ vertices from~$K$ to $x$.
Then 
\[
(i-j)X_{G_k}
+(j-k)X_{G_i}
+(k-i)X_{G_j}
=0,
\quad
\text{for any 
$0
\le i
\le j
\le k
\le \abs{K}$.}
\]
\end{corollary}
\begin{proof}
By \cref{prop:AP}, the chromatic symmetric functions $X_{G_j}$
form an arithmetic progression. In other words, if we let $m=\abs{K}$, 
then for any $0\le i\le m$, 
\begin{equation}\label{pf:Gijk}
X_{G_i}
=\frac{m-i}{m}X_{G_0}+\frac{i}{m}X_{G_m}.
\end{equation}
Let $F=(i-j)X_{G_k}
+(j-k)X_{G_i}
+(k-i)X_{G_j}$. Substituting \cref{pf:Gijk} and its analogues obtained 
by replacing $i$ with $j$ and $k$ respectively into $F$, we can
express $F$ as a linear combination of~$X_{G_0}$ and $X_{G_m}$.
In that expression, the coefficient of $X_{G_0}$ in $F$ is
\[
(i-j)\frac{m-k}{m}
+(j-k)\frac{m-i}{m}
+(k-i)\frac{m-j}{m}
=0, 
\]
and the coefficient of $X_{G_m}$ in $F$ is 
\[
(i-j)\frac{k}{m}
+(j-k)\frac{i}{m}
+(k-i)\frac{j}{m}
=0.
\]
This proves $F=0$ as desired.
\end{proof}

Here is an analogue of \cref{cor:AP:add}.
\begin{corollary}\label{cor:AP:remove}
Let $(K,x)$ be a tribe-vertex pair of a graph $H$.
Suppose that $x$ connects to all vertices of $K$.
Let~$H_i$ be the graph obtained from~$H$ by
removing $i$ edges from~$K$ to $x$.
Then
\[
(i-j)X_{H_k}
+(j-k)X_{H_i}
+(k-i)X_{H_j}
=0,
\quad\text{
for any 
$0
\le i
\le j
\le k
\le \abs{K}$.}
\]
\end{corollary}
\begin{proof}
Suppose that $\abs{K}=m$.
In \cref{cor:AP:add}, 
if we take $G=H$, then $H_i$ becomes $G_{m-i}$, and
\[
\brk1{(m-k)-(m-j)}X_{G_{m-i}}
+\brk1{(m-j)-(m-i)}X_{G_{m-k}}
+\brk1{(m-i)-(m-k)}X_{G_{m-j}}=0.
\]
Simplifying it yields the desired identity.
\end{proof}

\citet[Table~1]{SW16} discovered a captivating formula
using Stanley's generating function for Smirnov words, 
see also \citet[Theorem 7.2]{SW10}.
\begin{proposition}[\citeauthor{SW16}]\label{prop:path}
$X_{P_n}
=\sum_{I\vDash n}
w_I
e_I$
for any $n\ge 1$, where
\begin{equation}\creflabel[def]{def:w}
w_I
=i_1(i_2-1)(i_3-1)\dotsm(i_{-1}-1).
\end{equation}
\end{proposition}

\begin{figure}[htbp]
\begin{tikzpicture}[decoration=brace]
\draw (0,0) node {$K_m$};
\node[ball] (c) at (0: \r) {};
\node[ball] at (30: \r) {};
\node[ball] at (-30: \r) {};

\node (p1) at ($ (c) + (\r, 0) $) [ball] {};
\node[ellipsis] at ($ (p1) + (\r-\eps, 0) $) {};
\node[ellipsis] (e2) at ($ (p1) + (\r, 0) $) {};
\node[ellipsis] at ($ (p1) + (\r+\eps, 0) $) {};
\node (p2) at ($ (p1) + (2*\r, 0) $) [ball] {};
\draw [decorate] (4, -.2) -- (2, -.2);
\node[below=8pt, font=\footnotesize] at (e2) {$l$ vertices};

\draw[edge] (-75: \r) arc (-75: 75: \r);
\draw[edge] (c) -- ($(p1) + (\r-2*\eps, 0)$);
\draw[edge] (p2) -- ($(p2) - (\r-2*\eps, 0)$);

\begin{scope}[xshift=5cm]
\draw (0,0) node {$C_m$};
\node[ball] (c) at (0: \r) {};
\node[ball] at (30: \r) {};
\node[ball] at (-30: \r) {};

\node (p1) at ($ (c) + (\r, 0) $) [ball] {};
\node[ellipsis] at ($ (p1) + (\r-\eps, 0) $) {};
\node[ellipsis] (e2) at ($ (p1) + (\r, 0) $) {};
\node[ellipsis] at ($ (p1) + (\r+\eps, 0) $) {};
\node (p2) at ($ (p1) + (2*\r, 0) $) [ball] {};
\draw [decorate] (4, -.2) -- (2, -.2);
\node[below=8pt, font=\footnotesize] at (e2) {$l$ vertices};

\draw[edge] (-75: \r) arc (-75: 75: \r);
\draw[edge] (c) -- ($(p1) + (\r-2*\eps, 0)$);
\draw[edge] (p2) -- ($(p2) - (\r-2*\eps, 0)$);
\end{scope}

\begin{scope}[xshift=10cm]
\node (O) at (0,0) [ball]{};

\node (a1) at (\r*2, \r) [ball]{};
\node (b1) at (\r*2, 0) [ball]{};
\node (c1) at (\r*2, -\r) [ball]{};

\node (a2) at (\r*3, \r) [ellipsis]{};
\draw [decorate] (4, .8) -- (2, .8);
\node[below=8pt, font=\footnotesize] at (a2) {$a$ vertices};

\node (b2) at (\r*3, 0) [ellipsis]{};
\draw [decorate] (4, -.2) -- (2, -.2);
\node[below=8pt, font=\footnotesize] at (b2) {$b$ vertices};

\node (c2) at (\r*3, -\r) [ellipsis]{};
\draw [decorate] (4, -1.2) -- (2, -1.2);
\node[below=8pt, font=\footnotesize] at (c2) {$c$ vertices};

\node (a2l) at (\r*3-\eps, \r) [ellipsis]{};
\node (b2l) at (\r*3-\eps, 0) [ellipsis]{};
\node (c2l) at (\r*3-\eps, -\r) [ellipsis]{};

\node (a2r) at (\r*3+\eps, \r) [ellipsis]{};
\node (b2r) at (\r*3+\eps, 0) [ellipsis]{};
\node (c2r) at (\r*3+\eps, -\r) [ellipsis]{};

\node (a3) at (\r*4, \r) [ball]{};
\node (b3) at (\r*4, 0) [ball]{};
\node (c3) at (\r*4, -\r) [ball]{};

\draw[edge] (O) -- (a1) -- ($(a2l) + (-\eps, 0)$);
\draw[edge] (O) -- (b1) -- ($(b2l) + (-\eps, 0)$);
\draw[edge] (O) -- (c1) -- ($(c2l) + (-\eps, 0)$);

\draw[edge] ($(a2r) + (\eps, 0)$) -- (a3);
\draw[edge] ($(b2r) + (\eps, 0)$) -- (b3);
\draw[edge] ($(c2r) + (\eps, 0)$) -- (c3);
\end{scope}
\end{tikzpicture}
\caption{The lollipop $K_m^l$, 
the tadpole $C_m^l$,
and the $3$-spider $S(abc)$.}
\label{fig:lollipop-tadpole-3spider}
\end{figure}
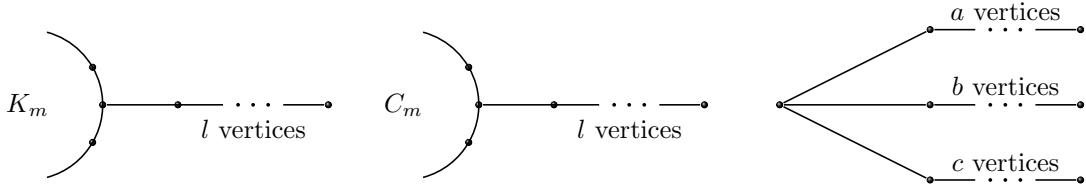

The \emph{lollipop} $K_m^l$ (resp., \emph{tadpole} $C_m^l$)
is the graph
obtained by identifying a vertex of the complete graph $K_m$ (resp., the cycle $C_m$)
and an end of the path $P_{l+1}$, 
see the leftmost (resp., middle) illustration in \cref{fig:lollipop-tadpole-3spider}.
\citet[Thm.~2]{Tom24} computed the chromatic symmetric function
of melting lollipops.
We restate his result for lollipops as a specialization. 

\begin{proposition}[Lollipops, \citeauthor{Tom24}]\label{prop:lollipop}
For any $n\ge a\ge 1$,
\[
X_{K_a^{n-a}}
=(a-1)!
\sum_{
\substack{
I\vDash n,\, i_{-1}\ge a
}}
w_I e_I.
\]
\end{proposition}

\citet{WZ25} provided a formula for tadpoles.
For any number $a\le n$, define
\begin{equation}\creflabel[def]{def:sigma}
\sigma_I(a)
=\min\{i_1+\dots+i_k\colon
0\le k\le \ell(I),\
i_1+\dots+i_k\ge a\}.
\end{equation}
The \emph{$a$-surplus} of $I$ is the number
\begin{equation}\creflabel[def]{def:Theta+}
\Theta_I(a)
=\sigma_I(a)-a.
\end{equation}

\begin{proposition}[Tadpoles, \citeauthor{WZ25}]\label{prop:tadpole}
For any $0\le l\le n-2$,
\[
X_{C_{n-l}^l}
=\sum_{I\vDash n}
\Theta_I(l+1)
w_I
e_I.
\]
\end{proposition}

For any composition $I\vDash n-1$,
the \emph{spider} $S(I)$ is the tree of order~$n$
obtained by identifying an end of the paths $P_{i_1+1}$, $P_{i_2+1}$, $\dots$.
The identified vertex is said to be the \emph{center} of $S(I)$.
We say that $S(I)$ is an $\ell(I)$-spider,
see the rightmost illustration in \cref{fig:lollipop-tadpole-3spider} 
for the $3$-spider $S(abc)$.
\citet[Lemma~4.4]{Zhe22} expressed
the chromatic symmetric function
of $3$-spiders in terms of chromatic symmetric functions of paths.

\begin{proposition}[$3$-Spiders, \citeauthor{Zhe22}]\label{prop:spider}
For any composition $abc\vDash n-1$, 
\[
X_{S(abc)}
=\sum_{i=0}^c
X_{P_i}
X_{P_{n-i}}
-\sum_{i=b+1}^{b+c}
X_{P_i}
X_{P_{n-i}}.
\]
\end{proposition}

Let $K_{i_1},\dots,K_{i_l}$ be completed graphs, where $i_k\ge 2$ for all $k$.
Suppose that $u_k$ and $v_k$ are distinct vertices of $K_{i_k}$.
The \emph{$K$-chain} $K_{i_1}+\dots+K_{i_l}$
is the graph obtained by identifying $v_k$ and $u_{k+1}$ for all $1\le k\le l-1$.
\citet[Thm.~3]{Tom24} discovered an elegant positive $e_I$-expansion for $K$-chains.
A \emph{weak composition} of~$n$
is a sequence of nonnegative integers whose 
sum is $n$.

\begin{theorem}[$K$-chains, \citeauthor{Tom24}]\label{thm:K-chain}
Let $\gamma=\gamma_1\dotsm\gamma_l$ be a composition with each part at least $2$.
If $K_\gamma$ is the $K$-chain with maximal clique sizes $\gamma_i$, then
\[
X_{K_\gamma}
=(\gamma_1-2)!\dotsm
(\gamma_{l-1}-2)!
(\gamma_l-1)!
\sum_{\alpha \in A_\gamma}
\brk4{
\alpha_1 
\prod_{i=2}^l
\abs{\alpha_i-\gamma_{i-1}+1}}
e_\alpha,
\]
where $l=\ell(\gamma)$ 
and $A_\gamma$ is the set of weak compositions 
$\alpha=\alpha_1\dotsm \alpha_l$ of $\abs{\gamma}-l+1$
such that for each $2\le i\le l$, either
\begin{itemize}
\item
$\alpha_i < \gamma_{i-1}$ and 
$\alpha_i+\dots+\alpha_l
<\gamma_i+\dots+\gamma_l-(l-i)$,
or
\item
$\alpha_i \ge \gamma_{i-1}$
and
$\alpha_i+\dots+\alpha_l
\ge\gamma_i+\dots+\gamma_l-(l-i)$.
\end{itemize}
\end{theorem}

\section{Path-conjoined graphs}\label{sec:path-conjoined}

\begin{figure}[h]
\begin{tikzpicture}
\node (c) at (0,0) [ellipsis]{};
\node (a1) at (-\r, 0) [ball]{};
\node (a2) at (-\r*2, 0) [ball]{};
\node (b1) at (\r, 0) [ball]{};
\node (b2) at (\r*2, 0) [ball]{};

\node[above right] at (a2) {$u$};
\node[above left] at (b2) {$v$};

\node (cr) at (0+\eps, 0) [ellipsis]{};
\node (cl) at (0-\eps, 0) [ellipsis]{};

\draw[edge] (a2) arc [start angle=0, end angle=360,
x radius=\r, y radius=\r*.75];
\draw[edge] (b2) arc [start angle=180, end angle=-180,
x radius=\r*.75, y radius=\r*.5];
\draw[edge] (a2) -- (a1) -- ($(cl) - (\eps, 0)$);
\draw[edge] (b2) -- (b1) -- ($(cr) + (\eps, 0)$);

\node at (-\r*3, 0) {$G$};
\node at (\r*2.75, 0) {$H$};
\node[below=4pt, font=\footnotesize] at (0, 0) {length $k$};
\end{tikzpicture}
\caption{The path-conjoined graph $P^k(G,H)$.}
\label{fig:GH}
\end{figure}
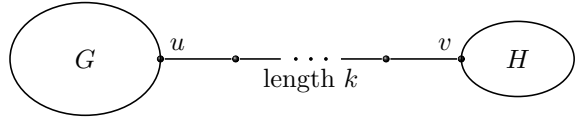

A \emph{rooted graph} $(G,u)$ is a pair 
consisting of a graph $G$ and a vertex $u$ of $G$.
We call $u$ the \emph{root} of~$G$.
For any rooted graphs $(G,u)$ and $(H,v)$, 
we define the \emph{path-conjoined graph} $P^k(G,H)$ 
to be the graph obtained by adding a path of length $k$
which links $u$ and $v$. It has order $\abs{G}+\abs{H}+k-1$,
see \cref{fig:GH}. 
We call $G$ and $H$ the \emph{node graphs} of $P^k(G,H)$.

A graph is \emph{vertex-transitive} if for any vertices~$x$ and~$y$,
there exists an automorphism $f$ such that $f(x)=y$.
We say that $G$ is a rooted graph without mentioning its root, 
whenever $G$ is vertex-transitive or 
the root is clear from context.
For instance,
path-conjoined graphs $P^k(K_g,K_h)$ are barbells,
see \cite[Theorem~3.7]{WZ25}. 
When one of the nodes, say, $H$, is a singleton $K_1$, we write
\[
G^k=P^k(G,K_1),
\]
and call it a \emph{tailed graph}.
Then lollipops are tailed cliques,
and tadpoles are tailed cycles.

The path-conjoining operation by a path of length $k$ 
is equivalent to joining a graph with a new edge $P_2$ for $k$ times,
as considered by \citet{GS01} and \citet{Tom24}.
As will be seen,
one benefits from the viewpoint of recognizing a path as a whole;
in this way the number of compositions
that appear in our formula decreases, 
in comparison with that in \cref{thm:K-chain}.

The path-conjoined graphs $P^k(G,H)$ with one node, say, $G$,
being a clique or a cycle is of particular interests.
We express the chromatic symmetric function of such a graph
in terms of those of tailed graphs with node $H$,
see \cref{prop:KPG,prop:CPG}.

%: KPG
\begin{proposition}[KPGs]\label{prop:KPG}
Let $H$ be a rooted graph.
Then for any $k\ge 0$ and $g\ge 1$,
\[
X_{P^k(K_g,\,H)}
=(g-1)!
\sum_{l=0}^{g-1}
(1-l)e_l
X_{H^{k+g-1-l}}.
\]
\end{proposition}
\begin{proof}
When $g=1$, 
the desired formula holds by the definition $G^k=P^k(G,K_1)$
for tailed graphs.
Let $g\ge 2$.
In \cref{cor:AP:remove},
if we take~$K=K_g$, 
and take $x$ to be the root of~$K_g$, 
then $G_0=G$, $G_{g-2}=P^{k+1}(K_{g-1},\,H)$, 
and $G_{g-1}$ is the disjoint union of~$K_{g-1}$ and~$H^k$.
It follows that 
\[
X_{P^k(K_g,\,H)}
=(g-1)
X_{P^{k+1}(K_{g-1},\,H)}
-(g-2)
X_{K_{g-1}}
X_{H^k}
\quad\text{for $g\ge 2$,}
\]
in which $X_{K_{g-1}}=(g-1)!e_{g-1}$.
The desired formula follows by using this relation iteratively.
\end{proof}

%: CPG
\begin{proposition}[CPGs]\label{prop:CPG}
Let $H$ be a rooted graph.
Then for any $k\ge 0$ and $g\ge 2$,
\[
X_{P^k(C_g,\,H)}
=(g-1)X_{H^{k+g-1}}
-\sum_{l=1}^{g-2}
X_{C_{g-l}}
X_{H^{k+l-1}}.
\]
\end{proposition}
\begin{proof}
The desired formula holds for $g=2$
since $H^{k+1}$ is the underlying graph 
of $P^k(C_2,\,H)$.
For $g\ge 3$,
let $u$ be the root of $C_g$, 
with neighbors $x$ and~$y$ on $C_g$.
Applying \cref{prop:3del} to the triangle $uxy$,
we obtain
\[
X_{P^k(C_g,\,H)}
=X_{P^{k+1}(C_{g-1},\,H)}
+X_{H^{k+g-1}}
-X_{C_{g-1}}
X_{H^k}
\quad\text{for $g\ge 3$}.
\]
We can then infer the desired formula by using this relation iteratively for $g$ down to $3$.
\end{proof}

The particular case for when $H=K_1$ appeared in 
\citet[Formula (3.11)]{LLWY21}.

A \emph{clique-path-cycle} graph is a graph of the form $P^b(K_a,C_c)$.
We call them KPC graphs for short.
Now we show their $e$-positivity 
by using \cref{prop:KPG} and the composition method.

%: KPC: c->c-1
\begin{theorem}[KPCs]\label{thm:KPC}
For any $a\ge 1$, $b\ge 0$ and $c\ge 2$,
the chromatic symmetric function
of the KPC graph $P^b(K_a,C_c)$ has the positive
$e_I$-expansion
\[
X_{P^b(K_a,\,C_c)}
=(a-1)!\sum_{K\vDash n}
c_K
w_K
e_K,
\]
where $n=a+b+c-1$, 
$w_K$ is defined by \cref{def:w}, and
\[
c_K=\begin{cases*}
0, 
& if $\ell(K)\ge 2$ and $k_2<a$,\\
\displaystyle
k_2-a-b
+\frac{k_2-k_1}{k_2-1},
& if $k_1\le a-1$ and $k_2\ge a+b$,\\
\Theta_K(a+b),
& otherwise.
\end{cases*}
\]
Here the function $\Theta_K(\cdot)$ is defined by \cref{def:Theta+}.
As a consequence, all KPC graphs are $e$-positive.
\end{theorem}
\begin{proof}
Let $G=P^b(K_a,\,C_c)$. 
Since $c\ge 2$,
by using \cref{prop:KPG,prop:tadpole},
we can deduce that
\[
\frac{X_G}{(a-1)!}
=\sum_{l=0}^{a-1}
(1-l)
e_l
X_{C_c^{a+b-1-l}}
=\sum_{K\vDash n}
\Theta_K(a+b)
w_K
e_K
-
\sum_{
K\vDash n,\
\ell(K)\ge 2,\
k_2\le a-1}
q_K
w_K
e_K,
\]
where $q_K=\Theta_{K\backslash k_2}(a+b-k_2)$.
Now we deal with the coefficient in the last sum.
\begin{enumerate}
\item
If $k_1<a+b-k_2$, then
$q_K=\Theta_K(a+b)$.
\item
If $a+b-k_2\le k_1<a+b$, then
$q_K
=k_1+k_2-(a+b)
=\Theta_K(a+b)$.
\item
If $k_1\ge a+b$, then 
$q_K
=k_1+k_2-(a+b)
=\Theta_K(a+b)+k_2$.
\end{enumerate}
Therefore,
the last sum can be rewritten as
\begin{align*}
\sum_{\substack{
K\vDash n,\
\ell(K)\ge 2\\
k_2\le a-1}}
q_K
w_K
e_K
&=\sum_{
\substack{
K\vDash n,\
\ell(K)\ge 2\\
k_2\le a-1}}
\Theta_K(a+b)
w_K
e_K
+\sum_{
\substack{
K\vDash n,\
\ell(K)\ge 2\\
k_1\ge a+b,\
k_2\le a-1}}
k_2
w_K
e_K.
\end{align*}
It follows that 
\begin{equation}\label{pf:KPC}
\frac{X_G}{(a-1)!}
=(c-1)ne_n
+\sum_{K\in\mathcal A_n}
\Theta_K(a+b)
w_K
e_K
-\sum_{K\in\mathcal B_n}
k_2
w_K
e_K,
\end{equation}
where
\begin{align*}
\mathcal A_n
&=\{K\vDash n\colon
\ell(K)\ge 2,\
k_2\ge a\},
\quad\text{and}\\
\mathcal B_n
&=\{K\vDash n\colon
\ell(K)\ge 2,\
k_1\ge a+b,\
2\le k_2\le a-1
\}.
\end{align*}
For any $K\in \mathcal B_n$, let $f(K)$ be the composition obtained 
by exchanging the first two parts of $K$.
Then~$f$ is a bijection between $\mathcal B_n$ and the set 
\[
f(\mathcal B_n)
=\{K\vDash n\colon \ell(K)\ge 2,\ 2\le k_1\le a-1,\ k_2\ge a+b\}.
\]
Then we can rewrite the negative sum in \cref{pf:KPC} as
\[
\sum_{
K\in \mathcal B_n}
k_2
w_K
e_K
=\sum_{K\in f(\mathcal B_n)}
\frac{k_1-1}
{k_2-1}
\cdotp
k_2
w_{K}
e_K.
\]
It is clear that $f(\mathcal B_n)\subseteq 
\mathcal A_n$.
For $K\in f(\mathcal B_n)$, we have
\[
\Theta_K(a+b)
-\frac{k_2(k_1-1)}{k_2-1}
=(k_2-a-b)
+\frac{k_2-k_1}{k_2-1}
>0.
\]
It follows that 
\begin{align*}
\frac{X_G}{(a-1)!}
&=(c-1)ne_n
+\sum_{K\in \mathcal A_n\backslash f(\mathcal B_n)}
\Theta_K(a+b)
w_K
e_K
+\sum_{K\in f(\mathcal B_n)}
\brk3{
k_2-a-b
+\frac{k_2-k_1}{k_2-1}}
w_K
e_K,
\end{align*}
which can be recast as the desired formula.
\end{proof}

\cref{thm:KPC} reduces to \cref{prop:lollipop} for $c=2$,
and to \cref{prop:tadpole} for both $a=1$ and $a=2$.
For example,
taking $(a,k,b)=(4,2,4)$, we can calculate
\begin{align*}
X_{P^2(K_4,C_4)}
&=18e_{1422}
+ 162e_{144}
+ 30e_{162}
+ 126e_{18}
+ 48e_{252}
+ 132e_{27}
+ 54e_{342}
+ 54e_{36}\\
&\qquad
+ 288e_{45}
+ 270e_{54}
+ 162e_9.
\end{align*}

A \emph{kayak paddle} is a 
path-conjoined graph with cycle nodes.
\citeauthor{AWv24} extended 
the appendability technique of $(e)$-positivity 
and proved the $e$-positivity of some chain-like graphs
including kayak paddles, 
see \cite[Proposition~34]{AWv24}.
\begin{proposition}[\citeauthor{AWv24}]\label{prop:CPC}
All kayak paddles are $e$-positive.
\end{proposition}

\citet{TW24X} discovered a positive $e_I$-expansion for the chromatic symmetric function of kayak paddles.

\section{Spider-conjoined graphs}\label{sec:spider-conjoined}
Let $\tau=\tau_1\dotsm\tau_l$ be a weak composition.
Let $\{s_1,\dots,s_l\}$ be the set of leaves of the spider~$S(\tau)$ with center $c$,
in which the leg containing~$s_i$ has length~$\tau_i$.
For any rooted graphs $(G_1,u_1),\dots,(G_l,u_l)$,
we define the $l$-\emph{spider-conjoined graph} $G=S^\tau(G_1,\dots,G_l)$ 
to be the graph obtained by identifying~$u_i$ with~$s_i$ for all $i$.
It has order 
\[
\abs{G}
=\abs{G_1}+\dots+\abs{G_l}
+\abs{\tau}
-l
+1.
\]
In particular, $1$-spider-conjoined graphs are tailed graphs,
and $2$-spider-conjoined graphs are path-conjoined graphs.
We say that $G$ has \emph{center} $c$ and \emph{nodes} $G_i$,
and omit the root of $G_i$ whenever $G_i$ is vertex-transitive
or the root is clear from context.

The chromatic symmetric function of 
$l$-spider-conjoined graphs for $l\ge 4$
can be expressed in terms of those for $3$-, $2$-, and $1$-spider-conjoined graphs,
when $\tau$ has at least two positive parts.

%: ell-spider-conjoined
\begin{proposition}\label{prop:ell-spider}
Let $\tau=\tau_1\dotsm \tau_l$ be a weak composition with $l\ge 3$ and $\tau_1,\tau_2\ge 1$.
Then for any $l$-spider-conjoined graph
$G=S^\tau(G_1,\dots,G_l)$,
\[
X_G
=X_{S^{(\tau_1-1)\tau_2 1}}
\brk1{
G_1,\,G_2,\,H}
+X_{G_1^{\tau_1-1}}
X_{S^{\tau_2\dotsm\tau_l}(G_2,\,\dots,\,G_l)}
-X_{P^{\tau_1+\tau_2-1}(G_1,G_2)}
X_H,
\]
where $H=S^{\tau_3\dotsm\tau_l}(G_3,\dots,G_l)$ is considered 
to be a rooted graph whose root is the center of $G$.
\end{proposition}
\begin{proof}
Applying \cref{fml:3del:12}
to the triangle formed by the center of $G$
and its neighbors on the $G_1$- and $G_2$-legs,
we obtain the desired formula.
\end{proof}

With \cref{prop:ell-spider} in hand,
one may concentrate on $3$-spider-conjoined graphs.

%: 3-spider-conjoined
\begin{proposition}\label{prop:3-spider}
Let $G$, $H$ and $J$ be rooted graphs.
Then for any $g,j\ge0$ and $h\ge 1$,
\[
X_{S^{ghj}(G,H,J)}
=X_{S^{(g+j)h0}(G,H,J)}
+\sum_{i=1}^j
\brk1{
X_{P^{g+h+i-1}(G,H)}
X_{J^{j-i}}
-
X_{G^{g+i-1}}
X_{P^{h+j-i}(H,J)}}.
\]
\end{proposition}
\begin{proof}
It is trivial for $j=0$. 
For $j\ge 1$,
taking $(G_1,G_2,G_3)=(J,H,G)$ in \cref{prop:ell-spider} yields
\[
X_{S^{ghj}(G,H,J)}
=X_{S^{(g+1)h(j-1)}(G,H,J)}
+X_{P^{g+h}(G,H)}
X_{J^{j-1}}
-X_{G^g}
X_{P^{h+j-1}(H,J)}.
\]
Using it iteratively for $j$ down to $1$, we can derive the desired formula.
\end{proof}

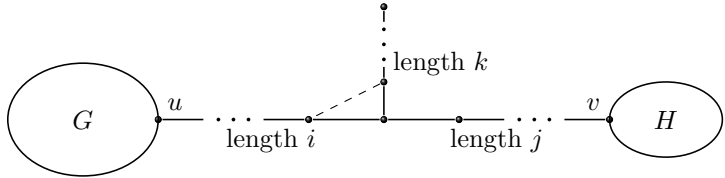
\begin{figure}[htbp]
\begin{tikzpicture}[decoration=brace]
\node (O) at (0,0) [ball]{};

\node (a1) at (-\r, 0) [ball]{};
\node (b1) at (\r, 0) [ball]{};
\node (c1) at (0, \r*.5) [ball]{};

\node (a2) at (-\r*2, 0) [ellipsis]{};
\draw [decorate] (0, -.2) -- (-3, -.2);
\node[below=8pt, font=\footnotesize] at (-\r*1.5, 0) {length $g$};

\node (b2) at (\r*2, 0) [ellipsis]{};
\draw [decorate] (3, -.2) -- (0, -.2);
\node[below=8pt, font=\footnotesize] at (\r*1.5, 0) {length $h$};

\node (c2) at (0, \r) [ellipsis]{};
\draw [decorate] (.2, 1.5) -- (.2, 0);
\node[right=8pt, font=\footnotesize] at (0, \r*.75) {length $j$};

\node (a2l) at (-\r*2+\eps, 0) [ellipsis]{};
\node (b2l) at (\r*2-\eps, 0) [ellipsis]{};
\node (c2b) at (0, \r-\eps) [ellipsis]{};

\node (a2r) at (-\r*2-\eps, 0) [ellipsis]{};
\node (b2r) at (\r*2+\eps, 0) [ellipsis]{};
\node (c2a) at (0, \r+\eps) [ellipsis]{};

\node (a3) at (-\r*3, 0) [ball]{};
\node[above right] at (a3) {$u$};
\node (b3) at (\r*3, 0) [ball]{};
\node[above left] at (b3) {$v$};
\node (c3) at (0, \r*1.5) [ball]{};

\draw[edge] (a3) arc [start angle=0, end angle=360,
x radius=\r, y radius=\r*.75];
\draw[edge] (b3) arc [start angle=180, end angle=-180,
x radius=\r*.75, y radius=\r*.5];
\draw[edge] (O) -- (a1) -- ($(a2l) - (-\eps, 0)$);
\draw[edge] (O) -- (b1) -- ($(b2l) + (-\eps, 0)$);
\draw[edge] (O) -- (c1) -- ($(c2b) + (0, -\eps*.5)$);;

\draw[edge] ($(a2r) - (\eps, 0)$) -- (a3);
\draw[edge] ($(b2r) + (\eps, 0)$) -- (b3);
\draw[edge] ($(c2a) + (0, \eps*.5)$) -- (c3);
\draw[dashed] (c1) -- (a1);

\node at (-\r*4, 0) {$G$};
\node at (\r*3.75, 0) {$H$};
\end{tikzpicture}
\caption{The spider-conjoined graph $S_j^{gh}(G,H)$.}
\label{fig:3spider-conjoined:GH}
\end{figure}

For any $1\le k\le l$, we denote 
\[
S^{\tau}(G_1,\dots,G_k,K_1,K_1,\dots,K_1)
=S^{\tau_1\dotsm\tau_k}_{\tau_{k+1}\dotsm\tau_l}(G_1,\dots,G_k)
\]
for compactness.
See \cref{fig:3spider-conjoined:GH} for the graph
$S^{gh}_j(G,H)$.

For when one or two of the nodes of a $3$-spider-conjoined graph
are the singletons $K_1$,
we write the specifications out.
It will be used in the proofs of \cref{thm:3spider:clique,thm:cycle2path:AB}.

%: 3-spider: K1 or K1K1
\begin{corollary}\label{cor:spider.ij:path}
Let $G$ and $H$ be rooted graphs.
Then for any $g,j\ge 0$ and $h\ge 1$,
\begin{align*}
X_{S_j^{gh}(G,H)}
&=\sum_{i=0}^j
X_{P_{i}}
X_{P^{g+h+j-i}(G,H)}
-\sum_{i=1}^{j}
X_{G^{g+i-1}}
X_{H^{h+j-i}},
\quad\text{and}\\
X_{S_{hj}^g(G)}
&=\sum_{i=0}^j
X_{P_i}
X_{G^{g+j-i+h}}
-\sum_{i=1}^{j}
X_{P_{i+h}}
X_{G^{g+j-i}}.
\end{align*}
\end{corollary}
\begin{proof}
They follow by taking $J=K_1$
and taking $J=H=K_1$ in \cref{prop:3-spider}, respectively.
\end{proof}

The second formula in \cref{cor:spider.ij:path} for $G=K_1$
reduces to \cref{prop:spider}.

In the remaining of \cref{sec:spider-conjoined},
we consider 
$3$-spider-conjoined graphs with a clique node 
and those with a cycle node;
see \cref{thm:3spider:clique,thm:3spider:cycle}, respectively.
The following \cref{lem:convolution:eP} will be used in proving 
\cref{thm:3spider:clique,thm:KmGK1}.
%and \cref{lem:sumsum:fgh} is for \cref{thm:3spider:clique}.

%: (1-l)e_lP_{n-l}
\begin{lemma}\label{lem:convolution:eP}
For $a\ge 0$,
\[
\sum_{l=0}^{a}
(1-l)e_l
X_{P_{n-l}}
=\begin{cases*}
\displaystyle
X_{K_{a+1}^{n-1-a}}\big/a!,
&\text{if $a\le n-1$};\\
e_n,
&\text{if $a=n$}.
\end{cases*}
\]
\end{lemma}
\begin{proof}
Let $0\le a\le n-1$.
By \cref{prop:path,prop:lollipop},
we can deduce that
\[
\sum_{l=2}^{a}
(l-1)e_l
X_{P_{n-l}}
=\smashoperator[r]{
\sum_{
K\vDash n,\ 
k_{-1}\le a}}
w_K
e_K
=X_{P_n}
-\smashoperator[r]{
\sum_{
K\vDash n,\ 
k_{-1}\ge a+1}}
w_K
e_K
=X_{P_n}
-X_{K_{a+1}^{n-a-1}}\big/a!.
\]
Rearranging the terms yields the desired formula for $a<n$.
Now, taking $a=n-1$, we obtain
\[
\sum_{l=0}^{n-1}
(1-l)e_l
X_{P_{n-l}}
=\frac{X_{K_n}}{(n-1)!}
=ne_n.
\]
Adding both sides by $(1-n)e_n$ yields the desired formula for $a=n$.
\end{proof}

Here is the reduction for the chromatic symmetric function of
$3$-spider-conjoined graphs $S^{ghk}(G,H,K)$
with a clique node $K=K_m$.

%: 3-spider: Km
\begin{theorem}\label{thm:3spider:clique}
Let $G$ and $H$ be rooted graphs.
For any $m, h\ge 1$ and $g, k\ge 0$, 
\begin{align*}
\frac{X_{S^{ghk}(G, H, K_m)}}{(m-1)!}
&=\sum_{z=0}^{m-1} e_{m-1-z} X_{P^{g+h+k+z}(G,H)}
-\sum_{\substack{a+b+c=k+m-2\\ a,\,b,\,c\ge 0,\ a\le m-1}}
(1-a)e_a X_{G^{b+g}} X_{H^{c+h}}\\
&\qquad
+\frac{1}{(m-1)!}\sum_{z=0}^{k-1} X_{K_m^{z}} X_{P^{g+h+k-z-1}(G,H)}.
\end{align*}
\end{theorem}
\begin{proof}
Applying \cref{prop:KPG} to the spider $S_0^{gh}(G,H)$
and the clique $K_m$, we obtain
\[
X_{S^{ghk}(G,\,H,\,K_m)}
=(m-1)!\sum_{l=0}^{m-1}
(1-l)
e_l
X_{S^{gh}_{k+m-1-l}(G,H)}.
\]
We write $q=k+m-1$ and $P^l(G,H)=P^l$ for short.
Using \cref{cor:spider.ij:path},
we can deduce that
\begin{equation}\label{pf:AB}
\frac{X_{S^{ghk}(G,\,H,\,K_m)}}{(m-1)!}
=\sum_{l=0}^{m-1}
(1-l)
e_l
\brk4{\sum_{z=0}^{q-l}
X_{P_{q-l-z}}
X_{P^{g+h+z}}
-\sum_{z=0}^{q-l-1}
X_{G^{g+q-l-z-1}}
X_{H^{h+z}}}.
\end{equation}
Let $A$ and $B$ be the first and second double sum on the right side of \cref{pf:AB}, respectively.
By interchanging the order of summation in $A$, we obtain
\begin{align*}
A&=\sum_{l=0}^{m-1} (1-l) e_l
\sum_{z=0}^{q-l} X_{P_{q-l-z}} X_{P^{g+h+z}}
=\sum_{z=0}^{q} X_{P^{g+h+z}}
\sum_{l=0}^{\min(m-1,\,q-z)} (1-l) e_l X_{P_{q-l-z}}
\\
&=\sum_{z=0}^{q} X_{P^{g+h+q-z}}
\sum_{l=0}^{\min(m-1,\,z)} (1-l) e_l X_{P_{z-l}}
\\
&=\sum_{z=0}^{m-1} X_{P^{g+h+q-z}}
\sum_{l=0}^{z} (1-l) e_l X_{P_{z-l}}
+\sum_{z=m}^{k+m-1} X_{P^{g+h+q-z}}
\sum_{l=0}^{m-1} (1-l) e_l X_{P_{z-l}}.
\end{align*}
Using \cref{lem:convolution:eP,prop:lollipop}, we can further transform it as 
\begin{align*}
A
&=\sum_{z=0}^{m-1}e_zX_{P^{g+h+q-z}}
+\sum_{z=m}^{k+m-1}X_{P^{g+h+q-z}}\frac{X_{K_m^{z-m}}}{(m-1)!}\\
&=\sum_{z=0}^{m-1}e_zX_{P^{g+h+q-z}}
+\frac{1}{(m-1)!}\sum_{z=0}^{k-1} X_{P^{g+h+q-z-m}} X_{K_m^{z}}.
%&=\sum_{z=0}^{m-1}e_zX_{P^{g+h+q-z}}
%+\sum_{z=0}^{k-1}\sum_{I\vDash m+z,\ i_{-1}\ge m}
%w_I e_I X_{P^{g+h+k-1-z}}\\
%&=\sum_{z=0}^{m-1}e_zX_{P^{g+h+k+m-1-z}}
%+\sum_{\abs{I}\le m+k-1,\ i_{-1}\ge m}
%w_I e_I X_{P^{g+h+k+m-1-\abs{I}}}.
\end{align*}
On the other hand, the double sum $B$ can be recast as the compact form
\[
B=\sum_{l=0}^{m-1}(1-l)e_l
\sum_{z=0}^{q-l-1}X_{G^{g+q-l-z-1}} X_{H^{h+z}}
=\sum_{\substack{a+b+c=k+m-2\\ a,\,b,\,c\ge 0,\ a\le m-1}}
(1-a)e_a X_{G^{b+g}} X_{H^{c+h}}.
\]
The sum $A+B$ is the desired formula.
\end{proof}

When $H=K_m$,
one may further simplify the formula in \cref{thm:3spider:clique}.

% 3spider: KmK1
\begin{corollary}\label{cor:3spider:GK}
Let $G$ be a rooted graph.
For any $k,g\ge 0$ and $m,h\ge 1$,
\begin{align*}
X_{S_h^{gk}(G,\,K_m)}
&=(m-1)!\brk4{
\sum_{z=0}^{m-1} e_z X_{G^{g+h+k+m-z-1}}
-\sum_{z=1}^{m-1} \frac{X_{K_z^h}}{(z-1)!} X_{G^{g+k+m-z-1}}}\\
&\qquad+
\sum_{z=0}^{k-1} 
\brk4{X_{K_m^{z}} X_{G^{g+h+k-z-1}}-X_{K_m^{h+z}} X_{G^{g+k-z-1}}}.
\end{align*}
\end{corollary}
\begin{proof}
In \cref{thm:3spider:clique}, taking $H=K_1$, we obtain
\begin{align*}
\frac{X_{S^{gk}_h(G, K_m)}}{(m-1)!}
&=\sum_{z=0}^{m-1}e_zX_{G^{g+h+k+m-1-z}}
-\sum_{\substack{a+b+c=k+m-2\\ a,\,b,\,c\ge 0,\ a\le m-1}}
(1-a)e_a X_{G^{b+g}} X_{P_{c+h+1}}\\
&\qquad
+\frac{1}{(m-1)!}\sum_{z=0}^{k-1} X_{K_m^{z}} X_{G^{g+h+k-z-1}}.\end{align*}
Let $A$ be the second sum, namely, the one with parameters $a$, $b$, and $c$. Then
\begin{align*}
A
&=\sum_{i=0}^{m-1}(1-i)e_i
\sum_{j=0}^{k+m-2-i}X_{G^{g+j}}X_{P_{h+k+m-i-j-1}}
\\
&=\sum_{j=0}^{k+m-2}
X_{G^{g+j}}
\sum_{i=0}^{\min(m-1,\,k+m-j-2)}(1-i)e_i X_{P_{h+k+m-i-j-1}}\\
&=\sum_{j=0}^{k+m-2}
X_{G^{g+k+m-j-2}}
\sum_{i=0}^{\min(m-1,\,j)}(1-i)e_i X_{P_{h-i+j+1}}.
\end{align*}
By \cref{lem:convolution:eP}, we can further derive that
\begin{align*}
A&=\sum_{j=0}^{k+m-2} X_{G^{g+k+m-j-2}}
\frac{X_{K_{\min(m,\,j+1)}^{h+j+1-\min(m,\,j+1)}}}{\min(m-1,\,j)!}\\
&=\sum_{j=0}^{m-2} X_{G^{g+k+m-j-2}} \frac{X_{K_{j+1}^h}}{j!}
+\sum_{j=m-1}^{k+m-2} X_{G^{g+k+m-j-2}} \frac{X_{K_m^{h+j-m+1}}}{(m-1)!}\\
&=\sum_{j=1}^{m-1} X_{G^{g+k+m-j-1}}\frac{X_{K_j^h}}{(j-1)!}
+\sum_{j=0}^{k-1} X_{G^{g+k-j+1}} \frac{X_{K_m^{h+j}}}{(m-1)!}.
\end{align*}
Substituting it back to the first formula in this proof, we obtain the desired one.
\end{proof}

The pineapple graph $\mathrm{Pi}_m^{\lambda}$ studied in \cite{TSH16}
is the graph $S^0_\lambda(K_m)$.
Now we can give a formula for pineapples with $\ell(\lambda)=2$. It is a graph obtained from a rooted clique $(K_m,u)$
by adding two disjoint paths with an end $u$,
whose lengths are $\lambda_1$ and $\lambda_2$ respectively.

\begin{corollary}\label{cor:pineapple:mab}
For any $g\ge 0$ and $h,m\ge 1$,
\[
X_{S_{gh}^0(K_m)}
=(m-1)!\brk4{\sum_{z=0}^{m-1} e_z X_{P_{g+h+m-z}}
-\sum_{z=1}^{m-1} \frac{X_{K_z^h}}{(z-1)!}X_{P_{g+m-z}}}.
\]
\end{corollary}
\begin{proof}
Immediately by taking $G=K_1$ and $k=0$ in \cref{cor:3spider:GK}.
\end{proof}

We end this section with a reduction for the chromatic symmetric function of
$3$-spider-conjoined graphs $S^{ghk}(G,H,K)$ 
with a cycle node $K=C_m$.

%: 3spider: Cm
\begin{theorem}\label{thm:3spider:cycle}
Let $G$ and $H$ be rooted graphs.
For any $g,h\ge 0$ and $m\ge 2$,
\[
X_{S^{ghk}(G,H,C_m)}
=(m-1)X_{S_{k+m-1}^{gh}(G,H)}
-\sum_{l=1}^{m-2}
X_{S_{k+l-1}^{gh}(G,H)}
X_{C_{m-l}}.
\]
\end{theorem}
\begin{proof}
In the graph $W=S_0^{gh}(G,H)$,
we set its root to be the spider center.
Then 
\[
W^k=S_k^{gh}(G,H)
\quad\text{and}\quad
P^k(W,C_m)=S^{ghk}(G,H,C_m).
\]
Taking $G=W$
in \cref{prop:CPG}
yields the desired formula.
\end{proof}

\section{Chain-conjoined graphs}\label{sec:chain-conjoined}

\begin{figure}[htbp]
\begin{tikzpicture}
\node[ball] (a1) at (\r*.75, 0) {};
\node[ball] (b1) at (2.75*\r, 0) {};
\node[ball] (b2) at (4.25*\r, 0) {};
\node[ball] (d1) at (8.25*\r+4*\eps, 0) {};
\node[ball] (c1) at (6.25*\r, 0) {};
\node[ellipsis] (c2) at (6.25*\r+\eps, 0) {};
\node[ellipsis] (c3) at (6.25*\r+2*\eps, 0) {};
\node[ellipsis] (c4) at (6.25*\r+3*\eps, 0) {};
\node[ball] (c5) at (6.25*\r+4*\eps, 0) {};
\draw[edge] (a1) arc [start angle=0, end angle=360,
x radius=\r*.75, y radius=\r*1.25];
\draw[edge] (b1) arc [start angle=180, end angle=-180,
x radius=\r*.75, y radius=\r*1.25];
\draw[edge] (d1) arc [start angle=180, end angle=-180,
x radius=\r*.75, y radius=\r*1.25];
\draw[edge] (a1) -- (b1);
\draw[edge] (b2) -- (c1);
\draw[edge] (d1) -- (c5);
\node[below=2pt, font=\footnotesize] at (1.75, 0) {length $\tau_1$};
\node[below=2pt, font=\footnotesize] at (5.25, 0) {length $\tau_2$};
\node[below=2pt, font=\footnotesize] at (7.25+4*\eps, 0) {length $\tau_l$};
\node at (0, 0) {$G_0$};
\node at (3.5*\r, 0) {$G_1$};
\node at (9*\r+4*\eps, 0) {$G_l$};
\end{tikzpicture}
\caption{The chain-conjoined graph $C(G_1, k_1, G_2, k_2, \dots ,G_{l-1}, k_{l-1}, G_l)$.}
\label{fig:chain-conjoined}
\end{figure}
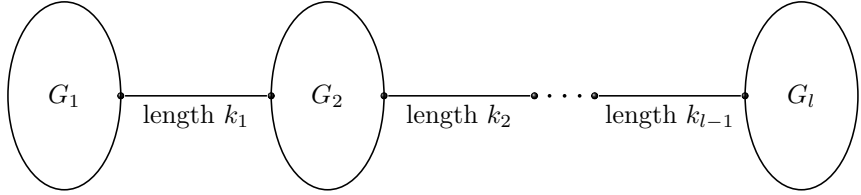

This section is devoted to chain-conjoined graphs.
A \emph{double rooted graph} $(G,u,v)$ is a triple
consisting of a graph $G$ and two vertices~$u$ and~$v$ of~$G$.
For our purpose, we restrict $u$ and $v$ to be distinct.
We call~$u$ and~$v$ its \emph{roots}.
For any double rooted graphs $(G_0,u_0,v_0)$, $\dots$, $(G_l,u_l,v_l)$,
and any weak composition $\tau=\tau_1\dotsm \tau_l$,
we define the $l$-\emph{chain-conjoined graph} 
$G
=
C^\tau
(G_0,\,
\dots,\,
G_l)
$
to be the graph obtained by adding a path 
of length~$\tau_i$ linking $v_i$ and $u_{i+1}$
for all $i$.
It has order 
\[
\abs{G}
=\abs{G_0}+\dots+\abs{G_l}+\abs{\tau}-l,
\]
see \cref{fig:chain-conjoined}. 
Then $1$-chain-conjoined graphs are path-conjoined graphs.
We call~$G_i$ the \emph{nodes} of $G$,
and omit the roots of $G_i$ 
whenever $G_i$ is vertex-transitive
or the roots of $G_i$ are clear from context.
Note that the roots $u_1$ and $v_l$ are useless
in our configuration.

We consider $2$-chain-conjoined graphs with the middle node $G_2=K_m$ for $m\ge 2$, denoted
\[
K_m^{gh}(G,H)
=C^{gh}(G,K_m,H).
\]

%: 3chain: Km
\begin{theorem}\label{thm:KGH}
Let $G$ and $H$ be rooted graphs.
Then for any $g,h\ge 0$ and $m\ge 2$,
\[
\frac{X_{K_m^{gh}(G,H)}}{(m-2)!}
=\sum_{z=1}^{m-1} z e_{m-1-z} X_{P^{g+h+z}(G,H)}
+(m-2)\sum_{\substack{a+b+c=m-2\\ a,\,b,\,c\ge 0}}
(a-1)e_a X_{G^{b+g}} X_{H^{c+h}}.
\]
\end{theorem}
\begin{proof}
The desired formula holds for $m=2$ since
$K_2^{gh}(G,H)
=P^{g+h+1}(G,H)$.
Suppose that $m\ge 3$.
We label the vertices of the clique $K_m$ in $G$ by
$uv_1v_2\dotsm v_{m-1}$, 
where~$u$ (resp., $v_{m-1}$)
is the root for conjoining~$K_m$ with~$G$ 
(resp., with~$H$).
Let $K_{m,k}^{gh}$ be the graph obtained from $K_m^{gh}=K_m^{gh}(G,H)$
by removing the edges $uv_1,\dots,uv_k$. 
By \cref{cor:AP:remove}, 
\[
X_{K_m^{gh}}
=(m-2)X_{K_{m,m-3}^{gh}}
-(m-3)X_{K_{m,m-2}^{gh}},
\quad\text{for $m\ge 3$}.
\]
Write $\mathfrak K_m^{gh}=S^{gh0}(G,H,K_m)$. Then the vertex in $K_m$ that connects $G$ with a path coincides the vertex that connects $H$ with a path, and
$K_{m,m-2}^{gh}=\mathfrak K_{m-1}^{(g+1)h}$.
Applying \cref{fml:3del:123} to the graph $K_{m,m-3}^{gh}$ with respect to the triangle $uvv_{m-2}$, we obtain
\[
X_{K_{m,m-3}^{gh}}
+X_{K_{m,\,m-1}^{gh}}
=X_{K_{m-1}^{(g+1)h}}
+X_{\mathfrak K_{m-1}^{(g+1)h}},
\quad\text{for $m\ge 3$}.
\]
Since $K_{m,\,m-1}^{gh}$
is the disjoint union of the graphs 
$G^g$ and $P^h(H,\,K_{m-1})$,
we can deduce that
\begin{multline*}
X_{K_m^{gh}}
=(m-2)X_{K_{m-1}^{(g+1)h}}
-(m-2)
X_{G^g}
X_{P^h(H,\,K_{m-1})}
+X_{\mathfrak K_{m-1}^{(g+1)h}}\\
=(m-2)
\brk2{
(m-3)X_{K_{m-2}^{(g+2)h}}
-(m-3)
X_{G^{g+1}}
X_{P^h(H,\,K_{m-2})}
+X_{\mathfrak K_{m-2}^{(g+2)h}}
}\\
-(m-2)
X_{G^g}
X_{P^h(H,\,K_{m-1})}
+X_{\mathfrak K_{m-1}^{(g+1)h}}
=\dotsm\\
=(m-2)!
X_{K_2^{(g+m-2)h}}
-\sum_{k=1}^{m-2}
(m-2)_k
X_{G^{g+k-1}}
X_{P^h(H,\,K_{m-k})}
+\sum_{k=1}^{m-2}
(m-2)_{k-1}
X_{\mathfrak K_{m-k}^{(g+k)h}}\\
=
\sum_{k=1}^{m-1}
(m-2)_{k-1}
X_{\mathfrak K_{m-k}^{(g+k)h}}
-
\sum_{k\ge 1}
(m-2)_k
X_{G^{g+k-1}}
X_{P^h(H,\,K_{m-k})}\\
=\sum_{k=0}^{m-2}
(m-2)_k
X_{S^{(g+k+1)h0}(G,H,K_{m-k-1})}
-\sum_{k\ge 1}
(m-2)_k
X_{G^{g+k-1}}
X_{P^h(H,\,K_{m-k})},
\end{multline*}
where $(m-2)_k=(m-2)(m-3)\dotsm(m-k-1)$ stands for the falling factorial.
By \cref{thm:3spider:clique,prop:KPG},
we can infer that
\begin{multline}\label{pf:Kmih}
X_{K_m^{gh}}
=(m-2)!
\sum_{k=0}^{m-2}
\sum_{z=0}^{m-k-2}
e_z
X_{P^{g+h+m-1-z}(G,H)}\\
-(m-2)!
\sum_{k=0}^{m-2}
\sum_{\substack{
a+b+c=g+h+m-2\\
a\ge 0,\,
b\ge g+k+1,\,
c\ge h}}
(1-a)
e_a
X_{G^b}
X_{H^c}\\
-(m-2)!
\sum_{k=1}^{m-2}
X_{G^{g+k-1}}
(m-k-1)
\sum_{l=0}^{m-k-1}
(1-l)e_l
X_{H^{h+m-k-1-l}}.
\end{multline}
Since the first double sum contains no parameter $k$, 
it can be simplified as
\[
\sum_{z=0}^{m-2}
(m-1-z)
e_z
X_{P^{g+h+m-1-z}(G,H)}
=\sum_{z=1}^{m-1}
z
e_{m-1-z}
X_{P^{g+h+z}(G,H)}.
\]
We write $Q=(1-a)
e_a
X_{G^b}
X_{H^c}$ for short.
Since the upper bound of $k$ for the second double sum is essentially $m-3$,
we can process the last two double sums by bootstrapping as 
\begin{align*}
&%1
\sum_{k=0}^{m-3}
\sum_{\substack{
a+b+c=g+h+m-2\\
a\ge 0,\,
b\ge g+k+1,\,
c\ge h}}
Q
+
\sum_{\substack{
a+b+c=g+h+m-2\\
a\ge 0,\,
b\ge g,\,
c\ge h
}}
(g+m-b-2)
Q\\
=&%2
\sum_{k=1}^{m-3}
\sum_{\substack{
a+b+c=g+h+m-2\\
a\ge 0,\,
b\ge g+k+1,\,
c\ge h}}
Q
+
\smashoperator[r]{
\sum_{\substack{
a+b+c=g+h+m-2\\
a\ge 0,\,
b\ge g+1,\,
c\ge h
}}}
(g+m-b-1)
Q
+(m-2)
\sum_{\substack{
a+c=h+m-2\\
a\ge 0,\,
c\ge h
}}
Q|_{b=g}\\
=&%3
\sum_{k=2}^{m-3}
\sum_{\substack{
a+b+c=g+h+m-2\\
a\ge 0,\,
b\ge g+k+1,\,
c\ge h}}
Q
+
\smashoperator[r]{
\sum_{\substack{
a+b+c=g+h+m-2\\
a\ge 0,\,
b\ge g+2,\,
c\ge h
}}}
(g+m-b)
Q
+(m-2)\brk4{
\smashoperator[r]{
\sum_{\substack{
a+c=h+m-3\\
a\ge 0,\,
c\ge h
}}}
Q|_{b=g+1}
+
\smashoperator[r]{
\sum_{\substack{
a+c=h+m-2\\
a\ge 0,\,
c\ge h
}}}
Q|_{b=g}
}\\
=&\dotsm\\
=&%5
\sum_{\substack{
a+b+c=g+h+m-2\\
a\ge 0,\,
b\ge g+m-2,\,
c\ge h
}}
(g+m-b+m-4)
Q
+(m-2)
\sum_{i=0}^{m-3}
\sum_{\substack{
a+c=h+m-2-i\\
a\ge 0,\,
c\ge h
}}
Q|_{b=g+i}
\\
=&%6
(m-2)
\sum_{i=0}^{m-2}
\sum_{\substack{
a+c=h+m-2-i\\
a\ge 0,\,
c\ge h
}}
Q|_{b=g+i}
=(m-2)
\sum_{k=0}^{m-2}
\sum_{a=0}^{m-2-k}
(1-a)
e_a
X_{G^{g+k}}
X_{H^{h+m-2-k-a}}.
\end{align*}
Substituting the previous two formulas into \cref{pf:Kmih},
we obtain the desired formula.
\end{proof}

\Cref{thm:KGH} reduces to a useful formula for PKPG
graphs by taking $H=K_1$.

\begin{theorem}[PKPGs]\label{thm:KmGK1}
For any $g,h\ge 0$ and $m\ge 2$,
\[
\frac{X_{K_m^{gh}(G,K_1)}}
{(m-2)!}
=\sum_{z=1}^{m-1}
z
e_{m-1-z}
X_{G^{g+h+z}}
-(m-2)
\sum_{z=1}^{m-1}
\frac{X_{K_z^h}}{(z-1)!} X_{G^{m-1-z+g}}.
\]
\end{theorem}
\begin{proof}
Taking $H=K_1$ in \cref{thm:KGH}, 
we obtain
\begin{align*}
\frac{X_{K_m^{gh}(G,K_1)}}{(m-2)!}
-\sum_{z=1}^{m-1}
z
e_{m-1-z}
X_{G^{g+h+z}}
&=(m-2)
\sum_{\substack{
a+b+c=g+h+m-2\\
a\ge 0,\,
b\ge g,\,
c\ge h
}}
(a-1)
e_a
X_{G^b}
X_{P_{c+1}}\\
&=(2-m)
\sum_{b=g}^{g+m-2}
X_{G^{b}}
\sum_{a=0}^{g+m-2-b}
(1-a)
e_a
X_{P_{g+h+m-1-a-b}}.
\end{align*}
Using \cref{lem:convolution:eP}, we can simplify the rightmost side of the equation above as
\[
(2-m)
\sum_{b=g}^{g+m-2}
X_{G^{b}}
\frac{X_{K_{g+m-2-b+1}^h}}{(g+m-2-b)!}
=(2-m)
\sum_{b=0}^{m-2}
\frac{X_{G^{m-2-b+g}}
X_{K_{b+1}^h}}
{b!}.
\]
Rearranging the terms yields the desired formula.
\end{proof}

By using the composition method with the aid of \cref{thm:KmGK1},
we now derive positive $e_I$-expansions for the chromatic symmetric functions 
of two kinds of $K$-chains,
see \cref{thm:PKP,thm:KKP}.
Let
\begin{equation}\creflabel[def]{def:f123}
\begin{cases}
f_1(I,a)=(a-1)w_I e_I,\\
f_2(I,a)=(a-2)i_{-1}w_{I\backslash i_{-1}}e_I,\quad\text{and}\\
f_3(I,a)=(i_{-1}-a+1)w_{I\backslash i_{-1}}e_I,
\end{cases}
\end{equation}
where the symbol $I\backslash i_{-1}$ stands for 
the composition obtained from $I$ by removing the last part.
Then 
\begin{equation}\label{f123}
f_1(I,a)
-f_2(I,a)
-f_3(I,a)
=\begin{dcases*}
(a-1)e_n,
& if $I=n$,\\
0,
& otherwise.
\end{dcases*}
\end{equation}
A \emph{path-clique-path} graph is a graph 
of the form $K_m^{gh}(K_1,K_1)$.

%: PKP
\begin{theorem}[PKPs]\label{thm:PKP}
Let $n=g+h+m$, where $g,h\ge 0$ and $m\ge 2$. 
Then the path-clique-path graph $K_m^{gh}(K_1,K_1)$ has the 
positive $e_I$-expansion
\[
\frac{X_{K_m^{gh}(K_1,K_1)}}
{(m-2)!}
=
(m-1)e_n
+
\sum_{
I\vDash n,\
\Theta_I(h+1)\ge m-1}
f_2(I,m)
+
\sum_{
I\vDash n,\
i_{-1}\ge m-1}
f_3(I,m),
\]
where $f_2$ and $f_3$ are defined by \cref{def:f123}.
\end{theorem}
\begin{proof}
Let $G=K_{m}^{gh}(K_1,K_1)$.
By \cref{thm:KmGK1,prop:path,prop:lollipop},
\begin{align*}
\frac{X_G}
{(m-2)!}
&=\sum_{z=1}^{m-1}
z
e_{m-1-z}
X_{P_{g+h+z+1}}
-(m-2)
\sum_{z=1}^{m-1}
\frac{X_{P_{m-z+g}}
X_{K_z^h}}{(z-1)!}\\
&=\sum_{z=0}^{m-1}
(m-1-z)
e_z
X_{P_{n-z}}
-(m-2)
\sum_{z=1}^{m-1}
\sum_{J\vDash m-z+g}
w_J
e_J
\smashoperator[r]{
\sum_{
I\vDash h+z,\
i_{-1}\ge z
}}
w_I
e_I.
\end{align*}
It follows that 
\begin{equation}\label{pf:PKP}
\frac{X_G}
{(m-2)!}
=\sum_{z=0}^{m-2}
\sum_{Iz\vDash n}
(m-1-z)
w_I
e_{Iz}
-(m-2)
\sum_{z=1}^{m-1}
\smashoperator[r]{
\sum_{
I\!J\vDash n,\
I\vDash h+z,\
i_{-1}\ge z
}}
w_I
w_J
e_{I\!J}.
\end{equation}
For the first double sum, we separate its summand for $z=0$
and rewrite it as
\begin{equation}\label{pf:PKP:1}
\sum_{z=0}^{m-2}
\sum_{Iz\vDash n}
(m-1-z)
w_I
e_{Iz}
=\sum_{K\vDash n}
(m-1)
w_K
e_K
+
\sum_{K\vDash n,\ 
k_{-1}\le m-2}
(m-1-k_{-1})
w_{K\backslash k_{-1}}
e_K.
\end{equation}
In the second double sum,
the composition $J$ can be equivalently replaced with its reversal $\overline{J}$; namely,
\begin{align*}
\sum_{z=1}^{m-1}
\sum_{
I\!J\vDash n,\
I\vDash h+z,\
i_{-1}\ge z}
w_I
w_J
e_{I\!J}
=\sum_{z=1}^{m-1}
\sum_{
I\!J\vDash n,\
I\vDash h+z,\
i_{-1}\ge z}
w_I
w_{\overline{J}}
e_{I\!J}.
\end{align*}
We regard $I\!J$ as a whole competition $K$.
On the one hand,
since $\abs{I}
=h+z
\le h+i_{-1}$,
we find $i_1+\dots+i_{-2}\le h$.
Since $\abs{I}=h+z\ge h+1$,
the composition~$I$ must be the shortest prefix of $K$
such that $\abs{I}\ge h+1$. 
On the other hand, 
$
w_I
w_{\overline{J}}
e_{I\!J}
=k_{-1}
w_{K\backslash  k_{-1}} 
e_K$.
As a result, we can rewrite the second double sum as
\begin{equation}\label{pf:PKP:2}
\sum_{z=1}^{m-1}
\sum_{
I\!J\vDash n,\
I\vDash h+z,\
i_{-1}\ge z
}
w_I
w_J
e_{I\!J}
=
\sum_{
K\vDash n,\
\Theta_K(h+1)\le m-2}
k_{-1}
w_{K\backslash  k_{-1}} 
e_K.
\end{equation}
Substituting \cref{pf:PKP:1,pf:PKP:2} into \cref{pf:PKP}, we obtain
\[
\frac{X_G}{(m-2)!}
=
\sum_{K\vDash n}
f_1(K,m)
-
\sum_{
K\vDash n,\
\Theta_K(h+1)\le m-2}
f_2(K,m)
-
\sum_{K\vDash n,\
k_{-1}\le m-2}
f_3(K,m).
\]
We proceed in $4$ cases. 
Let $K\vDash n$, $x=\Theta_K(h+1)$,
$y=k_{-1}$, and
$I=K\backslash k_{-1}$.
Denote by $c_K$ the coefficient of $e_K$ in $X_G/(m-2)!$.
Here we distinguish $c_K$ and $c_H$ if $K$ and $H$ are distinct compositions.
Then we have the following arguments.
\begin{enumerate}
\item
If $x,y\ge m-1$, then $c_K=(m-1)w_K$.
\item
If $x\ge m-1>y$, then $c_K=(m-1)w_K+(m-1-y)w_I=(m-2)yw_I$.
\item
If $x,y<m-1$, then $c_K=0$ by the previous item.
\item
If $x<m-1\le y$, then $c_K=(m-1)w_K-(m-2)yw_I=(y-m+1)w_I$.
\end{enumerate}
Gathering them together, we obtain
\[
\frac{X_{K_{m}^{gh}(K_1,K_1)}}
{(m-2)!}
=
\sum_{
\substack{
K\vDash n,\
k_{-1}\le m-2\\
\Theta_K(h+1)\ge m-1}}
f_2(K,m)
+
\sum_{
\substack{
K\vDash n,\ 
k_{-1}\ge m-1\\
\Theta_K(h+1)\ge m-1}}
f_1(K,m)
+
\sum_{
\substack{K\vDash n,\
k_{-1}\ge m-1\\
\Theta_K(h+1)\le m-2}}
f_3(K,m).
\]
By \cref{f123}, one may recast it as desired.
\end{proof}

For example, by \cref{thm:PKP} we can calculate
\[
X_{K_4^{21}(K_1,K_1)}
=24e_{142}
+16e_{151}
+24e_{16}
+2e_{421}
+6e_{43}
+48e_{52}
+30e_{61}
+42e_7.
\]

A \emph{clique-clique-path} graph is a graph 
of the form $K_{b+1}^{0a}(K_c,K_1)$.

%: KKP
\begin{theorem}[KKPs]\label{thm:KKP}
Let $n=a+b+c$, where $a\ge 0$, $b\ge 1$ and $c\ge 1$. 
Then the chromatic symmetric function of the 
clique-clique-path graph $K_{b+1}^{0a}(K_c,K_1)$
has the positive $e_I$-expansion
\[
\frac{X_{K_{b+1}^{0a}(K_c,K_1)}}{(b-1)!(c-1)!}
=
\sum_{
\substack{
K\vDash n\\
k_{-1}\ge b+c
}}
b
w_K
e_K
+
\sum_{
\substack{
K\vDash n,\
\ell(K)\ge 2,\
k_{-1}+k_{-2}\ge b+c\\
k_{-1}\le \min(b-1,\,c-1)
\ \text{or}\
k_{-1}\ge \max(b+1,\,c)
}}
\abs{b-k_{-1}}
w_{K\backslash k_{-1}}
e_K,
\]
where the symbol $I\backslash i_{-1}$ stands for 
the composition obtained from $I$ by removing the last part,
and the function $w$ is defined by \cref{def:w}.
\end{theorem}
\begin{proof}
Let $G=X_{K_{b+1}^{0a}(K_c,\,K_1)}$.
By \cref{thm:KmGK1,prop:lollipop},
\begin{align*}
\frac{X_G}
{(b-1)!}
&=\sum_{z=1}^{b}
z
e_{b-z}
X_{K_{c}^{a+z}}
-(b-1)
\sum_{z=1}^{b}
X_{K_{c}^{b-z}}
\sum_{
J\vDash a+z,\
j_{-1}\ge z
}
w_J
e_J\\
&=\sum_{z=1}^{b}
z
e_{b-z}
\sum_{I\vDash a+c+z,\
i_{-1}\ge c}
(c-1)!\,
w_I
e_I
-(b-1)(c-1)!
\sum_{z=1}^{b}
\sum_{\substack{
I\vDash b+c-z,\
i_{-1}\ge c
\\
J\vDash a+z,\
j_{-1}\ge z
}}
w_I
w_J
e_{I\!J}.
\end{align*}
Separating the term for $z=b$ from the first double sum,
we can rewrite it as
\[
\frac{X_G}
{(b-1)!(c-1)!}
=Y_1
+Y_2
-Y_3,
\]
where
$Y_1
=b\sum_{
K\vDash n,\
k_{-1}\ge c}
w_K
e_K$, 
\begin{align*}
Y_2
&=\sum_{z=1}^{b-1}
\sum_{\substack{
I\vDash a+c+z\\
i_{-1}\ge c}}
z
w_I
e_{i_1\dotsm i_{-2}(b-z)i_{-1}}
=
\sum_{\substack{
K\vDash n,\
k_{-1}\ge c\\ 
\ell(K)\ge 2,\
k_{-2}\le b-1}}
(b-k_{-2})
w_{K\backslash k_{-2}}
e_K,
\quad\text{and}\\ 
Y_3
&=(b-1)
\sum_{z=1}^b
\sum_{
\substack{
I\vDash a+b+1-z,\
i_{-1}\ge b-z+1\\
J\vDash c+z-1,\
j_{-1}\ge c
}}
w_I
w_J
e_{I\!J}.
\end{align*}
In the last expression, 
since $a+1\le\abs{I}\le a+i_{-1}$,
the composition~$I$ must be the shortest prefix of 
the composition $I\!J$
such that $\abs{I}\ge a+1$. It follows that 
\begin{equation}\label{pf:Y3:KKP}
Y_3=(b-1)
\sum_{\substack{
K=I\!J\vDash n,\
j_{-1}\ge c, \
\abs{I}=\sigma_K(a+1)
}}
w_I
w_J
e_K,
\end{equation}
where the function $\sigma_K(\cdot)$ is defined by \cref{def:sigma}.
We claim that
\begin{align}\label{pkkY3}
Y_3=
\sum_{
\substack{
K\vDash n,\
c\le k_{-1}\le b+c-1\\
\ell(K)\ge 2,\
k_{-1}+k_{-2}\ge b+c}}
(b-1)
k_{-1}
w_{K\backslash k_{-1}}
e_K
+
\sum_{
\substack{
K\vDash n,\
k_{-1}\ge c\\
k_{-1}+k_{-2}\le b+c-1}}
(b-1)
k_{-2}
w_{K\backslash k_{-2}}
e_K.
\end{align}
In fact, we separate the sum in \cref{pf:Y3:KKP} according to whether $\ell(J)=1$.
\begin{enumerate}
\item
For $\ell(J)=1$, we have
\[
w_I w_J=k_{-1}
w_{K\backslash k_{-1}},\quad
k_{-1}\le b+c-1,\quad
\ell(K)\ge2,
\quad\text{and}\quad
k_{-1}+k_{-2}\ge b+c.
\]
Conversely, if these three inequalities hold,
then $\abs{k_1\dotsm k_{-3}}\le a$.
Since $\abs{I}=\sigma_K(a+1)$,
we find $k_1\dotsm k_{-3}$ must be a proper prefix of $I$.
Since $J$ is not empty, we derive that $I=k_1\dotsm k_{-2}$ and $\ell(J)=1$.
Thus the terms in \cref{pf:Y3:KKP} with $\ell(J)=1$ sum to the first sum in \cref{pkkY3}.
\item
For $\ell(J)\ge 2$,
define $J'$ to be the composition 
obtained by moving $j_{-2}$ to the beginning of $J$, i.e.,
$J'=j_{-2}
j_1
j_2
\dotsm
j_{-3}
j_{-1}$.
When $J$ runs over the set $\{J\vDash n-\abs{I}\colon j_{-1}\ge c\}$,
so does $J'$, and vice versa.
If we write $K=IJ'$,
then 
\[
w_Iw_{J'}=k_{-2}
w_{K\backslash k_{-2}}
\quad\text{and}\quad
k_{-1}+k_{-2}\le b+c-1.
\]
Conversely, if $k_{-1}+k_{-2}\le b+c-1$,
then $\abs{k_1\dotsm k_{-3}}\ge a+1$,
which implies that $I$ is a prefix of $k_1\dotsm k_{-3}$, and $\ell(J)\ge 2$.
Thus the terms in \cref{pf:Y3:KKP} with $\ell(J)\ge 2$ sum to the second sum in \cref{pkkY3}.
\end{enumerate}
This proves \cref{pkkY3}.
Now we compute $Y_1-Y_3$. Since for $i\in\{1,2\}$,
\[
bw_K-(b-1)k_{-i}
w_{K\backslash k_{-i}}
=w_{K\backslash k_{-i}}
\brk1{
b(k_{-i}-1)
-(b-1)k_{-i}}
=(k_{-i}-b)
w_{K\backslash k_{-i}},
\]
we obtain $Y_1-Y_3=Q_1+Q_2+Q_3$, where
\begin{align*}
Q_1
&=\sum_{
K\vDash n,\
k_{-1}\ge b+c}
b
w_K
e_K,\\
Q_2
&=
\sum_{
\substack{
K\vDash n,\
c\le
k_{-1}\le b+c-1\\
\ell(K)\ge 2,\
k_{-1}+k_{-2}\ge b+c}}
(k_{-1}-b)
w_{K\backslash k_{-1}}
e_K,\quad\text{and}\\
Q_3
&=
\sum_{
\substack{
K\vDash n,\
k_{-1}\ge c\\
\ell(K)\ge 2,\
k_{-1}+k_{-2}\le b+c-1}}
(k_{-2}-b)
w_{K\backslash k_{-2}}
e_K.
\end{align*}
Note that in $Q_3$, we always have $k_{-2}\le b-1$.
Therefore,
\begin{multline*}
Y_2+Q_3
=
\sum_{\substack{
K\vDash n,\
k_{-1}\ge c\\ 
\ell(K)\ge 2,\
k_{-2}\le b-1}}
(b-k_{-2})
w_{K\backslash k_{-2}}
e_K
+\sum_{
\substack{
K\vDash n,\
k_{-1}\ge c\\
\ell(K)\ge 2,\
k_{-1}+k_{-2}\le b+c-1}}
(k_{-2}-b)
w_{K\backslash k_{-2}}
e_K\\
=
\sum_{\substack{
K\vDash n,\
k_{-1}\ge c\\ 
\ell(K)\ge 2,\
k_{-2}\le b-1\\
k_{-1}+k_{-2}\ge b+c}}
(b-k_{-2})
w_{K\backslash k_{-2}}
e_K
=
\sum_{\substack{
K\vDash n,\
\ell(K)\ge 2,\\
k_{-2}\le b-1\\
k_{-1}+k_{-2}\ge b+c}}
(b-k_{-2})
w_{K\backslash k_{-2}}
e_K.
\end{multline*}
We then deal with the summands in $Q_2$ with $k_{-1}\le b-1$
by exchanging $k_{-2}$ and $k_{-1}$. 
As a result, we can write
$Q_2=Q_{21}+Q_{22}$, where
\begin{align*}
Q_{21}
&=\sum_{
\substack{
K\vDash n,\
\max(b,c)\le
k_{-1}\le b+c-1\\
\ell(K)\ge 2,\
k_{-1}+k_{-2}\ge b+c}}
(k_{-1}-b)
w_{K\backslash k_{-1}}
e_K,
\quad\text{and}\\
Q_{22}
&=\sum_{
\substack{
K\vDash n,\ 
c\le k_{-2}\le b-1\\
\ell(K)\ge 2,\
k_{-1}+k_{-2}\ge b+c}}
(k_{-2}-b)
w_{K\backslash k_{-2}}
e_K.
\end{align*}
We observe that $Q_{22}$ cancels out a part of $Y_2+Q_3$ as
\begin{multline*}
(Y_2+Q_3)+Q_{22}
=
\sum_{\substack{
K\vDash n,\
\ell(K)\ge 2\\
k_{-2}\le b-1\\
k_{-1}+k_{-2}\ge b+c}}
(b-k_{-2})
w_{K\backslash k_{-2}}
e_K
+\sum_{
\substack{
K\vDash n,\ 
\ell(K)\ge 2\\
c\le k_{-2}\le b-1\\
k_{-1}+k_{-2}\ge b+c}}
(k_{-2}-b)
w_{K\backslash k_{-2}}
e_K\\
=\sum_{\substack{
K\vDash n,\
\ell(K)\ge 2\\
k_{-2}\le \min(b,c)-1\\
k_{-1}+k_{-2}\ge b+c}}
(b-k_{-2})
w_{K\backslash k_{-2}}
e_K
=\sum_{\substack{
K\vDash n,\
\ell(K)\ge 2\\
k_{-1}\le \min(b,c)-1\\
k_{-1}+k_{-2}\ge b+c}}
(b-k_{-1})
w_{K\backslash k_{-1}}
e_K.
\end{multline*}
The parameters $k_{-1}$ and $k_{-2}$ are exchanged again 
in the last step.
Collecting the results in
\[
\frac{X_G}
{(b-1)!(c-1)!}
=Q_1+Q_{21}+(Y_2+Q_3+Q_{22}),
\]
we obtain the desired formula.
\end{proof}

For example, by \cref{thm:KKP}, we can calculate
\[
X_{K_7^{01}(K_3,K_1)}
=48 e_{126}
+720 e_{162}
+1152 e_{171}
+1680 e_{18}
+192 e_{27}
+144 e_{36}
+1008 e_{72}
+1536 e_{81}
+2160 e_{9}.
\]

At the end of this paper,
we consider 
$2$-chain-conjoined graphs with the middle node $G_2=C_m$ for $m\ge 3$,
with an additional condition that the two nodes of $G_2$ are adjacent, denoted
\[
C_m^{gh}(G,H)
=C^{gh}(G,C_m,H).
\]
In particular,
the $2$-chain-conjoined graph $C_m^{ab}(K_1,K_1)$
is the hat graph $H_{a,m,b}$ studied in \cite{WZ25}.
The adjacency condition for the roots of $C_m$
is essential for the $e$-positivity.
For instance, the graph $G=C^{11}(K_1,C_4,K_1)$ in which
the roots of $C_4$ are non-adjacent is not $e$-positive as
\[
X_G
=18e_6+22e_{51}-2e_{42}+6e_{411}+9e_{33}+4e_{321}-2e_{222}+e_{2211}.
\]

%: GCH
\begin{theorem}[GCHs]\label{thm:cycle2path:AB}
Let $(C_m,x,y)$ be a double-rooted cycle with adjacent roots $x$ and $y$,
where $m\ge 2$.
For any rooted graphs $(G,v_G)$ and $(H,v_H)$
and any $g,h\ge 0$,
let $C_m^{gh}(G,H)$ be the graph obtained 
by adding a path of length $g$ linking $x$ and $v_G$
and another path of length $h$ linking $y$ and $v_H$.
Then
\[
X_{C_m^{gh}(G,H)}
=
\smashoperator[r]{
\sum_{\substack{
a+b=g+h+m-1\\
a\ge 0,\
b\ge g+h+1}}}
(a+1)
X_{P_a}
X_{P^b(G,H)}
+\smashoperator{
\sum_{\substack{
a+b+c=g+h+m-2\\
a\ge g,\
b\ge h,\
c\ge 2
}}}
X_{G^a}
X_{H^b}
X_{C_c}
-(m-2)
\smashoperator{
\sum_{\substack{
a+b=g+h+m-2\\
a\ge g,\
b\ge h
}}}
X_{G^a}
X_{H^b}.
\]
\end{theorem}
\begin{proof}
We omit the variables $G$ and $H$ for short and 
write 
\[
P^k(G,H)=P^k,\quad
S_m^{gh}(G,H)=S_m^{gh}
\quad\text{and}\quad
C_m^{gh}(G,H)
=C_m^{gh}.
\]
For convenience, 
we label the vertices of the clique $C_m$ in $G$ by
$v_0,\dots,v_{m-1}$, where $u=v_0$ and $v=v_{m-1}$.
Applying the triple-deletion property to the triangle $uvv_1$, we obtain
\[
X_{C_m^{gh}}
=X_{C_{m-1}^{(g+1)h}}
+X_{S_{m-2}^{(g+1)h}}
-X_{G^g}
X_{P^h(H,\,C_{m-1})},
\quad\text{for $m\ge 3$}.
\]
Iteratively using it, we can infer that
\begin{equation}\label{pf:GHC:0}
X_{C_m^{gh}}
=\sum_{k=1}^{m-1}
X_{S_{m-1-k}^{(g+k)h}}
-\sum_{k=1}^{m-2}
X_{G^{g+k-1}}
X_{P^h(H,\,C_{m-k})}.
\end{equation}
This leads us to deal with the graphs
$S_{m-1-k}^{(g+k)h}$ 
and $P^h(H,C_{m-k})$, respectively.
On one hand,
in \cref{cor:spider.ij:path}, 
replacing $k$ with $m-1-k$
and replacing $g$ with $g+k$,
we obtain
\[
X_{S_{m-1-k}^{(g+k)h}}
=\sum_{l=0}^{m-1-k}
X_{P_{m-1-k-l}}
X_{P^{g+h+k+l}}
-\sum_{l=0}^{m-2-k}
X_{G^{g+m-2-l}}
X_{H^{h+l}}.
\]
It follows that 
\begin{align}
\notag
\sum_{k=1}^{m-1}
X_{S_{m-1-k}^{(g+k)h}}
&=\sum_{k=1}^{m-1}
\sum_{l=k}^{m-1}
X_{P_{m-1-l}}
X_{P^{g+h+l}}
-\sum_{k=1}^{m-1}
\sum_{l=0}^{m-2-k}
X_{G^{g+m-2-l}}
X_{H^{h+l}}\\
\label{pf:GHC:1}
&=\sum_{l=1}^{m-1}
(m-l)
X_{P_{m-1-l}}
X_{P^{g+h+l}}
-\sum_{l=0}^{m-3}
(m-2-l)
X_{G^{g+m-2-l}}
X_{H^{h+l}}.
\end{align}
On the other hand, in \cref{prop:CPG},
taking $k=h$, $G=H$, and setting $m$ to be $m-k$, we obtain
\[
X_{P^h(H,\,C_{m-k})}
=(m-k-1)X_{H^{h+m-k-1}}
-\sum_{l=1}^{m-k-2}
X_{H^{h+l-1}}
X_{C_{m-k-l}}.
\]
It follows that 
\begin{align}
\notag
\smashoperator[r]{
\sum_{k=1}^{m-2}}
X_{G^{g+k-1}}
X_{P^h(H,C_{m-k})}
&=\smashoperator[r]{
\sum_{k=1}^{m-2}}
(m-k-1)
X_{G^{g+k-1}}
X_{H^{h+m-k-1}}
-
\smashoperator[l]{
\sum_{k=1}^{m-2}}
\smashoperator[r]{
\sum_{l=1}^{m-k-2}}
X_{G^{g+k-1}}
X_{H^{h+l-1}}
X_{C_{m-k-l}}\\
\label{pf:GHC:2}
&=\sum_{l=1}^{m-2}
l\cdotp
X_{G^{g+m-l-2}}
X_{H^{h+l}}
-\sum_{l=1}^{m-3}
\sum_{k=1}^{m-2-l}
X_{G^{g+k-1}}
X_{H^{h+l-1}}
X_{C_{m-k-l}}.
\end{align}
Substituting \cref{pf:GHC:1,pf:GHC:2} into \cref{pf:GHC:0}, we obtain
\begin{multline*}
X_{C_m^{gh}}
=\sum_{l=1}^{m-1}
(m-l)
X_{P_{m-1-l}}
X_{P^{g+h+l}}
-\sum_{l=0}^{m-3}
(m-2-l)
X_{G^{g+m-2-l}}
X_{H^{h+l}}\\
-\sum_{l=1}^{m-2}
l
X_{G^{g+m-l-2}}
X_{H^{h+l}}
+\sum_{l=1}^{m-3}
\sum_{k=1}^{m-2-l}
X_{G^{g+k-1}}
X_{H^{h+l-1}}
X_{C_{m-k-l}}\\
=\sum_{l=1}^{m-1}
(m-l)
X_{P_{m-1-l}}
X_{P^{g+h+l}}
+\smashoperator[l]{
\sum_{l=0}^{m-4}}
\smashoperator[r]{
\sum_{k=0}^{m-4-l}}
X_{G^{g+k}}
X_{H^{h+l}}
X_{C_{m-k-l-2}}
-(m-2)
\smashoperator{
\sum_{l=0}^{m-2}}
X_{G^{g+m-2-l}}
X_{H^{h+l}}\\
=\sum_{k=1}^{m-1}
k\cdotp
X_{P_{k-1}}
X_{P^{g+h+m-k}}
+\smashoperator{
\sum_{k+l\le m-4,\,k,l\ge0}}
X_{G^{g+k}}
X_{H^{h+l}}
X_{C_{m-k-l-2}}
-(m-2)
\sum_{l=1}^{m-1}
X_{G^{g+m-1-l}}
X_{H^{h+l-1}},
\end{multline*}
which simplies to the desired formula.
\end{proof}

A \emph{hat-chain} is a $(l+1)$-chain-conjoined graph of the form
$
C^{\tau_0\tau_1\dotsm\tau_l}(
K_1,\,
C_{m_1},\,
C_{m_2},\,
\dots,\,
C_{m_l},\,
K_1)$,
where $C_j$ are double-rooted cycles
with adjacent roots,
see \cref{fig:hat-chain}. 
\begin{figure}[htbp]
\begin{tikzpicture}[scale=.85, decoration=brace]
\node (0) at (-2,0) [ball]{};
\node[ellipsis] (e1) at (-1,0) {};
\node[ellipsis] (e1l) at ($ (e1) - (\eps, 0) $) {};
\node[ellipsis] (e1r) at ($ (e1) + (\eps, 0) $) {};
\node (1) at (0,0) [ball]{};
\node (m1) at (3,0) [ball]{};
\node[ellipsis] (a1) at (1.5,0) {};
\node[ellipsis] (a1l) at ($ (a1) - (\eps, 0) $) {};
\node[ellipsis] (a1r) at ($ (a1) + (\eps, 0) $) {};
\draw[edge] (1) arc [start angle=180, end angle=0, x radius=3*\r/2, y radius=.75*\r];
\node[ellipsis] (b1) at (4, 0) {};
\node[ellipsis] (b1l) at ($ (b1) - (\eps, 0) $) {};
\node[ellipsis] (b1r) at ($ (b1) + (\eps, 0) $) {};
\node[ball] (b2) at (5, 0) {};
\draw[edge] (b2) arc [start angle=180, end angle=0, x radius=3*\r/2, y radius=.75*\r];
\node[ellipsis] (c1) at (6.5,0) {};
\node[ellipsis] (c1l) at ($ (c1) - (\eps, 0) $) {};
\node[ellipsis] (c1r) at ($ (c1) + (\eps, 0) $) {};
\node[ball] (c2) at (8, 0) {};
\node[ellipsis] (c3) at (9.5,0) {};
\node[ellipsis] (c3l) at ($ (c3) - (\eps, 0) $) {};
\node[ellipsis] (c3r) at ($ (c3) + (\eps, 0) $) {};
\node[ball] (d1) at (11, 0) {};
\node[ellipsis] (d2) at (12.5, 0) {};
\node[ellipsis] (d2l) at ($ (d2) - (\eps, 0) $) {};
\node[ellipsis] (d2r) at ($ (d2) + (\eps, 0) $) {};
\node[ball] (d3) at (14, 0) {};
\node[ball] (f1) at (16, 0) {};
\node[ellipsis] (f2) at (15,0) {};
\node[ellipsis] (f2l) at ($ (f2) - (\eps, 0) $) {};
\node[ellipsis] (f2r) at ($ (f2) + (\eps, 0) $) {};
\draw[edge] (d1) arc [start angle=180, end angle=0, x radius=3*\r/2, y radius=.75*\r];

\draw[edge] (0) -- ($ (e1l) - (\eps, 0) $);
\draw[edge] ($ (e1r) + (\eps, 0) $) -- (1);
\draw[edge] (1) -- ($ (a1l) - (\eps, 0) $);
\draw[edge] ($ (a1r) + (\eps, 0) $) -- (m1) -- ($ (b1l) - (\eps, 0) $);
\draw[edge] ($ (b1r) + (\eps, 0) $) -- (b2);
\draw[edge] (b2) -- ($ (c1l) - (\eps, 0) $);
\draw[edge] ($ (c1r) + (\eps, 0) $) -- (c2) -- ($ (c3l) - (\eps, 0) $);
\draw[edge] ($ (c3r) + (\eps, 0) $) -- (d1) -- ($ (d2l) - (\eps, 0) $);
\draw[edge] ($ (d2r) + (\eps, 0) $) -- (d3);
\draw[edge] (d3) -- ($ (f2l) - (\eps, 0) $);
\draw[edge] ($ (f2r) + (\eps, 0) $) -- (f1);

\draw [decorate] (0, -.2) -- (-2, -.2);
\node[below=8pt, font=\footnotesize] at (-1, 0) {length $\tau_0$};
\draw [decorate] (3, -.2) -- (0, -.2);
\node[below=8pt, font=\footnotesize] at (1.5, 0) {length $m_1-1$};
\draw [decorate] (5, -.2) -- (3, -.2);
\node[below=8pt, font=\footnotesize] at (4, 0) {length $\tau_1$};
\draw [decorate] (8, -.2) -- (5, -.2);
\node[below=8pt, font=\footnotesize] at (6.5, 0) {length $m_2-1$};
\draw [decorate] (14, -.2) -- (11, -.2);
\node[below=8pt, font=\footnotesize] at (12.5, 0) {length $m_l-1$};
\draw [decorate] (16, -.2) -- (14, -.2);
\node[below=8pt, font=\footnotesize] at (15, 0) {length $\tau_l$};
\end{tikzpicture}
\caption{The hat-chain $C^{\tau_0\tau_1\dotsm\tau_l}(
K_1,\,
C_{m_1},\,
C_{m_2},\,
\dots,\,
C_{m_l},\,
K_1)$.}
\label{fig:hat-chain}
\end{figure}

\begin{conjecture}\label{conj:hat-chain}
All hat-chains are $e$-positive.
\end{conjecture}

\Citet{WZ25} confirmed its truth for $l=1$, i.e.,
when the hat-chain is a single hat.
When $l=2$ and $\tau_0=\tau_2=0$, 
the hat-chains becomes kayak paddles,
and \cref{conj:hat-chain} is true by \cref{prop:CPC}.
In the first version of our manuscript ({\tt arXiv:2406.01418v1}),
the conjecture for hat-chains was stated with the restriction $\tau_0=\tau_l=0$.
In fact, it can be considered
equivalent to \cref{conj:hat-chain} if one allows $m_1=m_l=2$.
Recently, \citet{TV25} found an involution proof 
for confirming \cref{conj:hat-chain}.

\section*{Acknowledgment}
The third author is appreciative of the hospitality of the LIGM
at Université Gustave Eiffel, 
where he was visiting Professor Jean-Yves Thibon in 2024.
This paper was written and submitted during that period.
The authors thank Watson Wang for catching an error in an early version of \cref{thm:3spider:clique}.

\bibliography{../csf.bib}
\end{document}